\documentclass{amsart}
\usepackage[utf8]{inputenc}
\usepackage{amsmath}
\usepackage{amsthm}
\usepackage[colorlinks=true, allcolors=blue]{hyperref}
\usepackage{cleveref}

\usepackage{MnSymbol}
\usepackage{wasysym}
\usepackage{graphicx}
\usepackage{tikz}
\usepackage{setspace}
\usepackage{float}

\newtheorem{ex}{Example}[section]

\author{Alfredo Roque Freire}
\author{Manuel A. Martins}
\title{Essay on modality across different logics} 

\newtheorem{definition}{Definition}[section]
\newtheorem*{definition*}{Definition}

\newcommand{\so}{\rightarrow}

\newcommand{\pair}[1]{\langle #1 \rangle}

\usepackage[backend=bibtex,style=alphabetic,maxbibnames=15,maxcitenames=6,dateabbrev=false]{biblatex}
\begin{document}

\maketitle

\begin{abstract}
In this paper, we deal with the problem of putting together modal worlds that operate in different logic systems.
When evaluating a modal sentence $\Box \varphi$, we argue that it is not sufficient to inspect the truth of $\varphi$ in accessed worlds (possibly in different logics). 
Instead, ways of transferring more subtle semantic information between logical systems must be established.
Thus, we will introduce modal structures that accommodate communication between logic systems by fixing a common lattice $L$ where different logics build their semantics. The semantics of each logic being considered in the modal structure is a sublattice of $L$.
In this system, necessity and possibility of a statement should not solely rely on the satisfaction relation in each world and the accessibility relation. The value of a formula $\Box \varphi$ will be defined in terms of a comparison between the values of $\varphi$ in accessible worlds and the common lattice $L$.
We will investigate natural instances where formulas $\varphi$ can be said to be necessary$/$possible even though all accessible world falsify $\varphi$. Finally, we will discuss frames that characterize dynamic relations between logic systems: classically increasing, classically decreasing and dialectic frames.
\end{abstract}

\section{Introduction}

Kripke semantics, originally proposed for Modal logic \cite{kripke2007semantical}, is very useful for dealing with several non-classical logics. In fact, most non-classical logics can be understood in a modal structure with classical worlds (e.g. Priest's approach in \cite{priest2008introduction}). Whats more, we can make use of the Kripke style semantic to study modality independently of the assumption of a specific background logic -- i.e. the logic regulating the non-modal connectives. For example, let us consider truth-values and propositional connectives operating in a fuzzy logic. If this fuzzyness is understood as uncertainty in reasoning, we can then use a Kripke frame to investigate modal properties of reasoning with uncertainty. 
Similar strategy can be used for other background logics, so one can produce a modal version for their chosen logic.

The study of these phenomena is often referred to as non-classical modal logic. Whether at some point modal operators where considered (e.g. by figures like Quine \cite[p. 80-94]{quine1986philosophy}) to be `deviant', they are now widely regarded as an integral part of the standard logical toolkit. In a recent edition dedicated to this subject \cite{igpl-special-edition}, Cintula and Weber assert that ``to be ‘non-classical’ in modal logic today means taking a much larger step beyond
the narrow confines of two-valued classical logic than it did when Quine wrote.'' But, even though the boarders of what is taken to be `deviant' have been restricted, some assumptions are still persistent. This is the case for the idea that (SLB) \textbf{a modal structure should operate in a single logic background}. Though we may consider varying backgrounds when dealing with different phenomena, one hardly investigates modal structures where worlds operating in a logic access worlds operating in a different logic. In this context, our first goal in this paper is to consider (G1) \textit{non-trivial modal structures where worlds can operate in different logics}.

Many logicians have indeed freed themselves from Quine's charge that (CLS) ``change of logic, change of subject'' \cite[p. 80]{quine1986philosophy} in a metatheoretical level, where one considers different non-classical modal logics for different purposes. But we have not so much freed ourselves from this charge when dealing with the inner aspects of a modal logic. In Berto and Jago recent book \cite{berto2019impossible} however, they explore ideas like \textit{impossible worlds} and the modal structure of \textit{imagination}, \textit{Omniscience} and \textit{relevance}.
The book extends and pays tribute to works of Kripke on non-normal worlds in \cite{KRIPKE2014206} and Priest (\cite{priest2016thinking, priest2016towards, priest1992}), Nolan (\cite{nolan1997impossible}), Zalta (\cite{zalta1997classically}), Rantala (\cite{rantala1982impossible}), among themselves and others on impossible worlds.
In this literature, worlds are devised into normal worlds which operate (usually) classically and non-normal (or impossible) worlds which operate in a non-standard way (see. \cite[p. 99-101]{berto2019impossible})\footnote{
Although semantics with impossible worlds consider what we may call ``worlds operating in different logics'', this is not their primary objective, nor they consider the question of modality across different logics in general. Kripke's non-normal frames actually operate classically with a change in the way modality is treated in the non-normal worlds (where everything is possible and nothing necessary). Moreover, Priest does not properly consider worlds operating in different logics -- every world operate paraconsistently, though normal worlds are those where no contradiction can be found. And in general, it is fair to say that this literature is concerned not with change in logic background, but with how modalities work when we consider accessibility not only to possible worlds, but also to worlds that one take to be impossible in some way.
}. In a more mathematical investigation, Martins et al in \cite{BMC14,MMB18} define equational hybrid logic, that adds hybrid machinary and equations (as atomic sentences) to modal syntax. 
Their semantics are structured in Kripke frames, but it is such that the worlds are endowed with particular algebras (see also \cite{Manzano2019} and \cite{DS07}). 
Although the logic is strictly the same in all worlds, as the algebra may vary from one world to another, the evaluation of equations depends on the algebra we have at the current world. 


Of course, an assumption like SLB may persist simply because it is true. But it can also have persisted for there is a reasonable \textit{fear of becoming trivial}. Though we do not subscribe to Quine's position in CLS, we still hold that `change of logic' should be a concern. And the question of what it means to say that `$\varphi$ is possible' or `$\varphi$ is necessary' should be investigated carefully. Consider we are in a world $w$ that accesses a world $w'$ \textbf{in a different logic} and $w'$ validates $\varphi$; should we conclude that `$\varphi$ is possible'? Indeed, the fact that $\varphi$ is valid in the accessed world should influence whether or not $\varphi$ is taken to be possible. But what is the extent of this influence? A second persistent assumption we shall investigate in this paper is that (NAW) \textbf{necessity means `true in all accessed worlds'} or the similar version (PAW) \textbf{possibility means `true in some accessed world'}.

Although in most cases the validity of $\varphi$ in $w'$ is sufficient to conclude that $\varphi$ is possible in a $w$ accessing $w'$, it can still be argued that there are more intricate cases where this connection cannot be established.
Addressing this issue is an important goal of this paper. So in the following sections we should produce (G2) \textit{natural cases in which NAW and PAW fails}. Although one may struggle with modal statements when worlds access worlds in different logics, the idea of accessing possible worlds in different logics becomes natural once one accepts the legitimacy of different logic backgrounds. So the issue is not so much the fact that one may allow accessibility to different logic worlds, but that it may be hard to make sense of it.
In this context, the truth value of a formula $\box \varphi$ will not be determined simply by the value of $\varphi$ concerning the logics in accessible worlds, but rather by relativising the value of $\varphi$ to the current world in a systematic and meaningful way.

Several reasons can be used to justify the development of modality in different logical backgrounds.
One may be interested in this phenomenon from a metaphysical perspective in case one subscribes to a pluralist view on logic; it can also be due to an interest in modal operations of imagination or kwoledge; or even as a new opportunity to understand connections between logical assumptions. This will not be our focus in this paper, as we elaborate primarily on \textbf{how} modality across different logics should operate. As a methodological desiderata, we should not consider the full generality of this phenomenon. For if one wish to reevaluate a widely held principle, one should find situations in which some better or refined understanding can be reached in holding a weaker or different principle. And this is what we propose. In a limited setting, where the assumptions on how the transfer of information between logics are fixed, we will develop meaningful cases where principles SLB, NAW and PAW fail. 
Finally, we will briefly discuss the notion of frames appropriate to the proposed semantics, with examples of classic increasing, decreasing and dialectic frames.

\section{Necessity as true in all possible worlds}

Whenever logicians teach modal systems, they may start exploring students' intuitions on the contrast between `true but it may be false' and `true and it cannot be false'.
Stronger rhetorical force is achieved with standard examples such as ``is it or not necessarily true that creatures with a heart and those with kidneys are the same creatures''? Or the question explored in Frege's seminal work \textit{{\"U}ber sinn und bedeutung} \cite{frege1892sinn}: how can the sentence ``\textit{Hesperus} is the same planet as \textit{Phosphorus}'' be informative if both names refer to Venus?
It is later instigating the limits of our imagination that we may teach what it takes for something to be not only true, but necessarily so: `Can you imagine a scenario in which $\varphi$ is false? Does it mean that all imaginable scenarios are such that $\varphi$ is true?'
Students are quick to learn that necessity means `true in all possible worlds'. 

Consider we are dealing here with the idea of possible world as `those one can imagine'.
An important intuition that implicitly allows the student to conclude NAW is that worlds operate in the same intensional vocabulary. Though `having heart' and `having kidney' designate the same objects in our world, they express different senses (in Frege's \cite{frege1892sinn}) or intensions (in Carnap's \cite{Carnap.1947})\footnote{Barcan in \cite{barcan1946functional} extends this approach to first order modal theories. For a more recent approach, see \cite{fitting2004first}.}. And these intensions are preserved across possible worlds -- e.g. `having a heart' cannot mean `having a wing' in a possible world. 
One may reject a possible world in which the intension of a `triangle' is a `square' as it is used to evaluate the sentence `it is possible that a triangle has four sides'. 
In this understanding, imagination should be bounded, restricted. Otherwise, the analysis of modal statements trivializes -- i.e. all statements are possible and none necessary. On the issue of preserving meaning of basic vocabulary, Kripke comments \cite[p. 77]{Kripke1972}:

\begin{quote}
[\ldots] I don't mean, of course, that there mightn't be counterfactual situations in which in the other possible worlds people actually spoke a different language. One doesn't say that `two plus two equals four' is contingent because people might have spoken a language in which `two plus two equals four' meant that seven is even.
\end{quote}

The fear of becoming trivial when basic meanings are not fixed relates to an important intuition we have about modalities. Consider we say `imagine a world in which horses have big and strong wings, this being possible to imagine results in the possibility that horses can fly'. At least two questions arise when dealing with this statement: (i) can we still say that these beings are horses? (ii) if possibility is so free from restrains, then what explanatory role would it have?  If we want to produce meaningful analysis of modal statements and accommodate basic intuitions on canonical examples, 
the meaning of some terms (quantifiers, relations, names, functions) should be fixed or partially fixed as we navigate through possible worlds.

Ask the question `can you imagine a scenario in which $\varphi$ is false?', but now stressing that possibility extends to changes in logical background. If the assumed background logic is classical and $\lnot \varphi$ is provable in it, then one may answer negatively to the possibility of $\varphi$ for they reject the very idea of imagining a world in a different logic. But even if one accept that imagination can range over worlds in a different logics, the answer to the possibility of $\varphi$ may not be simple. 
One may accept imagining a world in a different logic while still being skeptical about how this imagination \textbf{connects to modal statements}. For this prudent character, one cannot say that $\varphi$ is possible \textbf{simply} by imagining a world where $\varphi$ is valid but in a different logic. As a residual implication of Quine's charge that `change of logic is change of subject', one should carefully consider what the statement $\varphi$ means in this other logic from ones point of view.

Both in the traditional modal logic as in the case we have worlds in different logics modalities \textit{make sense} or produce \textit{good explanations} so long as \textit{something relevant} is preserved across worlds. Even though Newton is the (or one of the) inventor of Calculus, the phrase (1) `It is possible that Newton died in infancy' might be judge as true while (2) `It is possible that the inventor of calculus died in infancy' is judge as false. This occurs, as suggested by Kripke (\cite{Kripke1972}) and Barcan (\cite{barcan1946functional}) theories, for names designate rigidly. This is a regularity that allows for an important and intuitive distinction between phrases (1) and (2). It is in virtue of preservation (more broadly or narrowly) of names and properties across worlds that the traditional modal analysis become meaningful.
And we assume it is for this reason that the subject of modality across varying logics is hardly considered. It is indeed reasonable that one understands that a predicate means different things in different logic backgrounds, and so one jump (prematurely, as we shall see) to the conclusion that logic cannot meaningfully vary between possible worlds.


But let us now consider a more deeply rooted assumption, one that forces us to impose all these regularities: we characterize necessity as being true in every accessible world. 
In \cite{Kripke1972}, Kripke considers the possibility that `Nixon could have a different name' (say Nixton). In this context, if this is a world in which Nixton lost the election, one may conclude (from the current world) that this world validates the statement `Nixon could have lost the election'. Kripke explains \cite[p. 49]{Kripke1972}:
\begin{quote}
[\ldots] proper names are rigid designators, for although the man (Nixon) might not have been the President, it is not the case that he might not have been Nixon (though he might not have been called `Nixon')
\end{quote}
Kripke assumes here that there are different levels of analysis for the name `Nixon', one that rigidly designate and another which is given by facts of the world like the `name used to call this person'. And the way he avoids getting lost in knowing which is the appropriate rigid name is saying that possible worlds are stipulated by the relevant changes we will find in such world (e.g. `the world in which Nixon is called Nixton and he lost the election').

However, we can alternatively say that there is a function across worlds making the appropriate translation of the `names used to call the person'\footnote{Fitting uses a similar approach to define his intensional logic for first order in \cite{fitting2004first}.}. In this case, the sentence `Nixton lost the election' validates `it is possible that Nixon lost the election' because sentences in a possible world are \textit{understood through a transference of semantic information} not necessarily originated in the evaluation of the same formula. And the stipulation of the alternative world is that `it is a world with a person called Nixton that lost the election' together with the stipulation that `Nixton is mapped to the person called Nixon in the current world'\footnote{
    Note that once we accept this flexibility, the question of the adequacy of stipulations arises. One might argue that there is no need to consider stipulations for the map `Nixon $\so$ Nixton' as the map will not be correct otherwise. This flexibility, therefore, makes Kripke's rigid designation thesis an explicit theoretical choice not previously incorporated into the modal framework. This is similar to a shift from `necessity means truth in all worlds' to `necessity means truth in all \textbf{accessed} worlds' even if we accept that possible worlds are only those accessed. In this case, one is just breaking the metaphysical position into more nuanced and separate theses.
}. What is at play here is not so much the preservation of names, but the agreement on a reference vocabulary and how information is transferred from one world to another.

For most familiar cases, fixing names is sufficient to provide a satisfactory account of modal distinctions of phrases involving names (even if there are alternative approaches). We can however investigate cases where the transfer of information between worlds is more complicated. 
Consider we inspect a sentence $\psi$ in a possible world and, instead of having a situation where a strategy for matching the vocabularies (e.g. `Nixton $\so$ Nixon') is available, there is no adequate translation that fully preserve the intended meaning of $\psi$ in the original world. 
This is often the case where the background logic of worlds are different. How can we preserve meaning of a formula if the basic logical operations do not work the same way, how can we translate statements of infinitary logic to traditional classical logic? 
But, though the meaning of a sentence may not be completely preserved, one can find ways to bring information of what occurs in the other world to their familiar understanding.

Before discussing worlds with different logics, let us consider the following meta-metaphysical structure that exemplifies our phenomena.
Take $w_a$ to be a world that accept the existence of abstract objects only when there are concrete objects to abstract from, and let $w_n$ be a world accessible to $w_a$ that rejects the existence of abstracts. If $w_n$ is a world where there are beings that are rational and winged, one may assert that $w_n$ validates $\Diamond \delta \equiv$`` the possibility of the existence of an abstract winged rational being'' in $w_a$ even though $\delta$ is false in $w_n$. 
Naturally, one may say that, if there is a world like $w_n$ in which there are rational winged beings, then there is also a world $w_n'$ that is identical to $w_n$ but in which one can abstract from collections of concretes.
Of course, the existence of $w_n'$ may avoid the need to inspect the `problematic' $w_n$, but it does not answer to the question ``does $w_n$ validate $\Diamond \delta$?'' Accepting that $w_n$ validate $\Diamond \delta$ in $w_a$ is intuitive however (or at least not counter-intuitive). One is in this case simply accepting that inspecting the possibility of an abstract $P$ can occur by observing concrete individuals that bear the property $P$ in a possible world. 

Now we consider a case where possible worlds operate in different logics. Let $w_p$ paraconsistent world accessible to a classical world $w_c$ in which a formula $\varphi$ has value $\{T, F\}$. How should the classical world $w_c$ understand $w_p$'s evaluation of $\varphi$ with respect to the possibility of $\varphi$? If $w_c$ takes $w_p$'s answer at face value, they would both evaluate that $w_p$ validates $\varphi$ is possible and does not validate that $\varphi$ is possible -- which, \textbf{one may argue}, does not ``make sense'' for the classical world. Instead, $w_c$ \textbf{can} understand $\{T, F\}$ as a \textbf{parameter to be accommodated in $\mathbf{w_c}$'s own semantics}. This can be done, for example (though not as we will propose next sections), by saying that the \{T, F\} being impossible in `my point of view' means that the statement is understood as `false in an accessed world'. Much more is required though to argue for a particular way of using parametric values in a different logic. And we shall develop a natural way of making this connection. The relevant aspect here is that one should find a way to or agree on the way information is transferred between worlds in a model before evaluating sentences in possible worlds. 

Necessity should be relativized to the semantics of the world evaluating the modal sentence. It is in this opening that NAW may fail. Even if an accessed world evaluates a formula as false/true, its evaluation when brought to the current world can render a different meaning than `false/true in an accessed world'.
Now, as we move to this understanding of necessity, we should ask: can we actually find a situation in which it make sense to say that $\varphi$ is necessary even though it is false in a accessed world? This question is relevant for, even though we do not assume NAW, there might not be natural cases in which NAW fails. We shall however develop intuitive models in which this occurs in the following sections.  

\section{Many-logic modal structure over lattices}

Many propositional logic systems have lattice-based semantics. This is the case of Classical Propositional Logic, which has Boolean Algebras as its equivalent algebraic semantics. A similar phenomenon happens with Intuitionistic Propositional Logic and Modal Logic which have, respectively, Heyting Algebras and Modal Algebras as their equivalent algebraic semantics (see \cite{BP89}). In all these cases the associated algebras are lattice expansions (i.e. lattices with extra operations).

Other systems do not have as strong a connection with lattices as Algebrizable Logics, while still having lattice semantics.
For example, the semantics of Logic of Paradox can be accounted with the gap and glut lattice; the Three Valued Logic has the three linear lattice as semantics; the $\vee\wedge$-fragment of Classical Propositional Logic, that does not have theorems, is associated with the variety of Distributive Lattices (\cite{Font1991}). 
\L ukasiewicz Logic is another example for which  the lattice $([0,1], \leq)$ is the basis of its semantics.
All of these examples use the meet (product ``.'') and join (sum ``+'') operation to represent conjunction and disjunction. If the value of $\phi$ is $x$ and of $\psi$ is $y$, then the value of $\phi \land \psi$ is $x.y$ and of $\phi \lor \psi$ is $x+y$.
Other logical operations are, in this case, defined according to the particular purposes of each logic. In the end, some values of the lattice are taken as those that determine whether a formula is considered valid or invalid. We call the set of these values a filter (in line with Blok and Pigozzi's nomenclature in \cite{BP89}) so that $\phi$ is valid ($\vDash \phi$) when the value assigned to $\phi$ in a valuation is in the filter.

Of course there are unorthodox logics for which no lattice can be used as semantics. An example of this phenomena can be obtained with any equational logic induced by a class of algebras which are not related with lattices (e.g. the equational logic of the variety of groups).
However, if we intend to develop a general framework for a modal logic in which we have different logics in different worlds, a common lattice together with an appropriate filter can unify the analysis of modal sentences. Different sublattices of a common lattice share a sense of order and joint/meet operations. 
So, in our system, information will be transferred between logic systems through the common order of the base lattice.
For operational reasons, we assume that lattices associated with logics are complete\footnote{A lattice $L$ is complete when all subsets of $L$ have meet and join. This condition trivially holds if $L$ is finite.}, and all of them are sublattices of a fixed lattice not necessarily complete.\footnote{Note that any two logics with lattice semantics have at least one lattice that have those as sublattices -- namely, the lattice that associate a value above all maximal values and a value bellow all minimal values of the lattices representing each logic.} The complete sublattices of the basis lattice represent logics in which our many-logic-modal universe will operate.

Although this choice simplifies our definitions, it is sufficiently generic to capture many logical systems in a natural and meaningful way. 

\begin{ex}
Consider the following lattice $L$
with at least the following three sublattices. Each of those sublattices are related to a logic system:

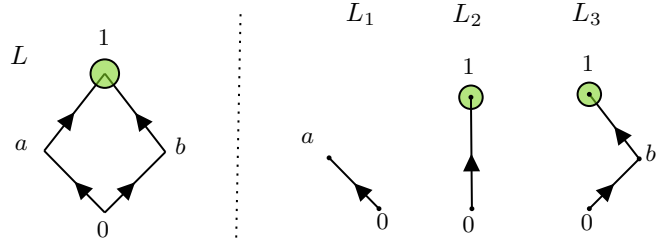
\begin{figure}[H]
\begin{center}
\tikzset{every picture/.style={line width=0.75pt}} 
\begin{tikzpicture}[x=0.75pt,y=0.75pt,yscale=-1,xscale=1]

	\draw    (164.86,144.22) -- (134.3,113.06) ;
	\draw [shift={(149.58,128.64)}, rotate = 45.56] [fill={rgb, 255:red, 0; green, 0; blue, 0 }  ][line width=0.08]  [draw opacity=0] (8.93,-4.29) -- (0,0) -- (8.93,4.29) -- cycle    ;
	\draw    (164.86,144.22) -- (195.43,113.66) ;
	\draw [shift={(180.14,128.94)}, rotate = 135] [fill={rgb, 255:red, 0; green, 0; blue, 0 }  ][line width=0.08]  [draw opacity=0] (8.93,-4.29) -- (0,0) -- (8.93,4.29) -- cycle    ;
	\draw    (195.43,113.66) -- (164.86,74.1) ;
	\draw [shift={(180.14,93.88)}, rotate = 52.31] [fill={rgb, 255:red, 0; green, 0; blue, 0 }  ][line width=0.08]  [draw opacity=0] (8.93,-4.29) -- (0,0) -- (8.93,4.29) -- cycle    ;
	\draw    (134.3,113.06) -- (164.86,74.1) ;
	\draw [shift={(149.58,93.58)}, rotate = 128.12] [fill={rgb, 255:red, 0; green, 0; blue, 0 }  ][line width=0.08]  [draw opacity=0] (8.93,-4.29) -- (0,0) -- (8.93,4.29) -- cycle    ;
	\draw  [fill={rgb, 255:red, 126; green, 211; blue, 33 }  ,fill opacity=0.58 ] (157.67,74.1) .. controls (157.67,70.13) and (160.89,66.91) .. (164.86,66.91) .. controls (168.83,66.91) and (172.05,70.13) .. (172.05,74.1) .. controls (172.05,78.08) and (168.83,81.3) .. (164.86,81.3) .. controls (160.89,81.3) and (157.67,78.08) .. (157.67,74.1) -- cycle ;
	\draw    (303.14,142.27) -- (277.97,116.61) ;
	\draw [shift={(290.56,129.44)}, rotate = 45.56] [fill={rgb, 255:red, 0; green, 0; blue, 0 }  ][line width=0.08]  [draw opacity=0] (8.93,-4.29) -- (0,0) -- (8.93,4.29) -- cycle    ;
	\draw  [fill={rgb, 255:red, 126; green, 211; blue, 33 }  ,fill opacity=0.58 ] (343.6,86.02) .. controls (343.6,82.75) and (346.25,80.1) .. (349.52,80.1) .. controls (352.79,80.1) and (355.44,82.75) .. (355.44,86.02) .. controls (355.44,89.29) and (352.79,91.94) .. (349.52,91.94) .. controls (346.25,91.94) and (343.6,89.29) .. (343.6,86.02) -- cycle ;
	\draw    (409.22,142.27) -- (434.38,117.11) ;
	\draw [shift={(421.8,129.69)}, rotate = 135] [fill={rgb, 255:red, 0; green, 0; blue, 0 }  ][line width=0.08]  [draw opacity=0] (8.93,-4.29) -- (0,0) -- (8.93,4.29) -- cycle    ;
	\draw    (434.38,117.11) -- (409.22,84.54) ;
	\draw [shift={(421.8,100.82)}, rotate = 52.31] [fill={rgb, 255:red, 0; green, 0; blue, 0 }  ][line width=0.08]  [draw opacity=0] (8.93,-4.29) -- (0,0) -- (8.93,4.29) -- cycle    ;
	\draw  [fill={rgb, 255:red, 126; green, 211; blue, 33 }  ,fill opacity=0.58 ] (403.3,84.54) .. controls (403.3,81.27) and (405.95,78.62) .. (409.22,78.62) .. controls (412.49,78.62) and (415.14,81.27) .. (415.14,84.54) .. controls (415.14,87.81) and (412.49,90.46) .. (409.22,90.46) .. controls (405.95,90.46) and (403.3,87.81) .. (403.3,84.54) -- cycle ;
	\draw    (350.01,143.01) -- (349.52,86.02) ;
	\draw [shift={(349.76,114.52)}, rotate = 89.5] [fill={rgb, 255:red, 0; green, 0; blue, 0 }  ][line width=0.08]  [draw opacity=0] (8.93,-4.29) -- (0,0) -- (8.93,4.29) -- cycle    ;
	\draw  [fill={rgb, 255:red, 0; green, 0; blue, 0 }  ,fill opacity=1 ] (277.24,116.61) .. controls (277.24,116.21) and (277.57,115.88) .. (277.97,115.88) .. controls (278.38,115.88) and (278.7,116.21) .. (278.7,116.61) .. controls (278.7,117.02) and (278.38,117.34) .. (277.97,117.34) .. controls (277.57,117.34) and (277.24,117.02) .. (277.24,116.61) -- cycle ;
	\draw  [fill={rgb, 255:red, 0; green, 0; blue, 0 }  ,fill opacity=1 ] (348.79,86.02) .. controls (348.79,85.62) and (349.11,85.29) .. (349.52,85.29) .. controls (349.92,85.29) and (350.25,85.62) .. (350.25,86.02) .. controls (350.25,86.43) and (349.92,86.75) .. (349.52,86.75) .. controls (349.11,86.75) and (348.79,86.43) .. (348.79,86.02) -- cycle ;
	\draw  [fill={rgb, 255:red, 0; green, 0; blue, 0 }  ,fill opacity=1 ] (302.41,142.27) .. controls (302.41,141.87) and (302.73,141.54) .. (303.14,141.54) .. controls (303.54,141.54) and (303.87,141.87) .. (303.87,142.27) .. controls (303.87,142.67) and (303.54,143) .. (303.14,143) .. controls (302.73,143) and (302.41,142.67) .. (302.41,142.27) -- cycle ;
	\draw  [fill={rgb, 255:red, 0; green, 0; blue, 0 }  ,fill opacity=1 ] (349.28,142.28) .. controls (349.28,141.88) and (349.61,141.55) .. (350.01,141.55) .. controls (350.41,141.55) and (350.74,141.88) .. (350.74,142.28) .. controls (350.74,142.68) and (350.41,143.01) .. (350.01,143.01) .. controls (349.61,143.01) and (349.28,142.68) .. (349.28,142.28) -- cycle ;
	\draw  [fill={rgb, 255:red, 0; green, 0; blue, 0 }  ,fill opacity=1 ] (408.49,142.27) .. controls (408.49,141.87) and (408.82,141.54) .. (409.22,141.54) .. controls (409.62,141.54) and (409.95,141.87) .. (409.95,142.27) .. controls (409.95,142.67) and (409.62,143) .. (409.22,143) .. controls (408.82,143) and (408.49,142.67) .. (408.49,142.27) -- cycle ;
	\draw  [fill={rgb, 255:red, 0; green, 0; blue, 0 }  ,fill opacity=1 ] (433.65,117.11) .. controls (433.65,116.7) and (433.98,116.38) .. (434.38,116.38) .. controls (434.79,116.38) and (435.11,116.7) .. (435.11,117.11) .. controls (435.11,117.51) and (434.79,117.84) .. (434.38,117.84) .. controls (433.98,117.84) and (433.65,117.51) .. (433.65,117.11) -- cycle ;
	\draw  [fill={rgb, 255:red, 0; green, 0; blue, 0 }  ,fill opacity=1 ] (408.49,84.54) .. controls (408.49,84.14) and (408.82,83.81) .. (409.22,83.81) .. controls (409.62,83.81) and (409.95,84.14) .. (409.95,84.54) .. controls (409.95,84.95) and (409.62,85.27) .. (409.22,85.27) .. controls (408.82,85.27) and (408.49,84.95) .. (408.49,84.54) -- cycle ;
	\draw  [dash pattern={on 0.84pt off 2.51pt}]  (233.32,44.1) -- (231.08,159.74) ;

	\draw (117.91,106.11) node [anchor=north west][inner sep=0.75pt]  [font=\small]  {$a$};
	\draw (198.81,104.91) node [anchor=north west][inner sep=0.75pt]  [font=\small]  {$b$};
	\draw (159.26,147.46) node [anchor=north west][inner sep=0.75pt]  [font=\small]  {$0$};
	\draw (159.86,50.38) node [anchor=north west][inner sep=0.75pt]  [font=\small]  {$1$};
	\draw (262.37,102.6) node [anchor=north west][inner sep=0.75pt]  [font=\small]  {$a$};
	\draw (300.34,143.06) node [anchor=north west][inner sep=0.75pt]  [font=\small]  {$0$};
	\draw (343.76,145.03) node [anchor=north west][inner sep=0.75pt]  [font=\small]  {$0$};
	\draw (343.76,65.1) node [anchor=north west][inner sep=0.75pt]  [font=\small]  {$1$};
	\draw (436.03,108.52) node [anchor=north west][inner sep=0.75pt]  [font=\small]  {$b$};
	\draw (403.46,143.55) node [anchor=north west][inner sep=0.75pt]  [font=\small]  {$0$};
	\draw (403.95,63.62) node [anchor=north west][inner sep=0.75pt]  [font=\small]  {$1$};
	\draw (284.99,37.14) node [anchor=north west][inner sep=0.75pt]  [font=\normalsize]  {$L_{1}$};
	\draw (338.77,37.64) node [anchor=north west][inner sep=0.75pt]  [font=\normalsize]  {$L_{2}$};
	\draw (398.97,37.14) node [anchor=north west][inner sep=0.75pt]  [font=\normalsize]  {$L_{3}$};
	\draw (115.63,58.78) node [anchor=north west][inner sep=0.75pt]  [font=\normalsize]  {$L$};
\end{tikzpicture}
\end{center}
\caption{Four valued base lattice.}
\label{four-valued-base-lattice}
\end{figure}
\end{ex}

In this example, the logic obtained from the values $\{0, a, 1\}$ is precisely the same as that obtained from $L_3$. 
However, their relationships to the values of $L_1$ or $\{a,1\}$ are different.
This results from the way we understand the values of $L_3$ in the bigger picture of $L$. 
Selecting the appropriate sublattice for a logic in this framework is not just the job of representing logical properties, but also of positioning the logical values in the common order of $L$.
Classical logic may be the logic of $L_2$, but it is also the logic of $\{a,1\}$ and $\{b, 1\}$. However, considering which is the correct choice for representing \textbf{the} Classical Logic in $L$ is not the purpose of this article. 
One could even argue that all three options are Classical Logic, adding that a particular choice of one of these options represents an additional metalogical and/or metaphysical commitment.
It suffices for our purpose to consider all complete sublattices as representing a logic with no specific commitment to actually representing \textbf{the} Classical Logic, or \textbf{the} Logic of Paradox, etc.

We restrict our investigation to a propositional language containing only conjunction, disjunction and negation. The base lattice $L$ will thus have an unary operation that we call ``complement" (represented by the minus symbol `$-$'). If a formula $\varphi$ has value $x \in L$, then $\lnot \varphi$ will have value $-x$.
With this addition, we should carefully consider how negation propagates to sublattices of $L$. One way is to require that our sublattices have the extra property of being closed under complements -- i.e. if $L'$ is a sublattice of $L$ and $x$ is a value in $L'$, then $-x$ is also a value in $L'$. In this case, we say that the negation is \textbf{rigidly interpreted}.

However, sometimes it is worth not having that extra condition on the sublattices (e.g. having more variety of sublattices/logics). In such case, we must be able to define an interpretation of the negation in each sublattice, induced by the negation of the basis lattice. A given sublattice $L'$ of $L$ may be such that, though $x \in L'$, it does not have $-x$. For this, we use the common order of $L$ to introduce our \textbf{down/up-interpretation} for values not present in a sublattice:

\begin{definition}
In a base lattice $L$, a value $a \in L$ is interpreted in a sublattice $L'$ as:
\begin{enumerate}
    \item \textbf{Down-interpretation} - the least value in $L'$ that is larger than \textbf{all values of $L'$ that are smaller than $a$} -- formally, $a^{L'} = \bigcup_{L'} \{x \in L' \mid x \leq a\}$ (if $\{x \in L' \mid x \leq a\} =\emptyset$, then $a^{L'}$ is the least value in $L'$).
    \item \textbf{Up-interpretation} - the largest value in $L'$ that is smaller than \textbf{all values of $L'$ that are bigger than $a$} -- formally, $a^{L'} = \bigcap_{L'} \{x \in L' \mid x \geq a\}$ (if $\{x \in L' \mid x \geq a\} =\emptyset$, then $a^{L'}$ is the greatest value in $L'$).
\end{enumerate}
\end{definition}

When the value $a$ is already in $L'$, then $a^{L'}$ is precisely $a$ in both the down and up interpretations. Also, while down-interpretation tends to make the value smaller if it is not present in the lattice, up-interpretation tends to make the value larger. We will assume the down-interpretation for the remainder of this article, as it can be understood as `more conservative'.

The interpretation of negation can now be stated precisely even when $x \in L'$ and $-x \notin L'$:\footnote{In what follows, if nothing is said, negation is rigidly interpreted.}
$$(-a)^{L'} = \bigcup_{L_1} \{x\in L_1 \mid x \leq \neg a\}$$

With the environment for logics to coexist and exchange information, we define our modal structures:

\begin{definition}
For a given lattice $L$ and propositional language $\mathcal{L}$, a \textbf{many-logic modal} $L$-structure $M$ is the tuple $\pair{W, R, s}$ such that
\begin{enumerate}
    \item each world $w_i$ in $W$ is the pair $\pair{i, L_i}$ where $i$ is the identifier for the world and $L_i$ is the complete sublattice of $L$ associated with the world $w_i$.
    \item $R$ is a relation from worlds to worlds, i.e. $R \subseteq W \times W$.
    \item and $s$ is a function from propositional variables to values in $L_i$ for each world $w = \pair{i, L_i}$, i.e. 
    $s: W \times Var \longrightarrow L$ such that $s(w_i, x) \in L_i$.
\end{enumerate}
\end{definition}

With $M = \pair{W,R,s}$ a many-logic modal $L$-structure, we define the valuation $v_{w}$ of propositional formulas in $\mathcal{L}$ in the world $w \in W$:

\newcounter{count-def-valuation}
\begin{definition}\label{valuation}
The valuation function $v_{w_i}: Form \longrightarrow L_i$ is such that 
\begin{enumerate}
    \item $v_{w_i}(x) = s(w_i, x)$.
    \item $v_{w_i}(\lnot \varphi) = (- v_{w_i}(\varphi))^{L_i}$.
    \item $v_{w_i}(\varphi \lor \psi) = (v_{w_i}(\varphi) + v_{w_i}(\psi))^{L_i}$.
    \item $v_{w_i}(\varphi \land \psi) = (v_{w_i}(\varphi) . v_{w_i}(\psi))^{L_i}$.
    \setcounter{count-def-valuation}{\value{enumi}}
\end{enumerate}

\end{definition}

The value of a formula $\Box \varphi$ in a world $w$ depends on the value $\varphi$ in the accessible worlds. Also, the idea of `all worlds' in the traditional definition equates in our system to the `meet operation'. Thence, we evaluate $\Box \varphi$ by taking the meet operation over the set of relativized values of $\varphi$ in all accessed worlds. Formally we have:

\begin{definition*}[Continuation of \cref{valuation}]\
\begin{enumerate}
    \setcounter{enumi}{\value{count-def-valuation}}
    \item 
    
    $v_{w_i}(\Box \varphi) = \bigcap\limits_{L_{w_i}}\{(v_{w}(\varphi))^{L_i} \mid w_i R w \land w \in W\} 
    $. Or, expanding with definition of down-interpretation:
    $$v_{w_i}(\Box \varphi) =\bigcap_{L_{w_i}} \{\bigcup_{L_{w_i}}\{x\in L_{w_i} \mid x \leq v_{w}(\varphi) \}:w_i R w \land w \in W\}$$
    \item $v_{w_i}(\Diamond \varphi) = (- v_{w_i}(\Box \lnot \varphi))^{L_i}$.\footnote{We consider $\Box$ equivalent to $\lnot \Diamond \lnot$ for simplicity. It may be defined, for instance, using directly up or down interpretations.}
\end{enumerate}
\end{definition*}

Let us evaluate an example of modal formulas:

\begin{ex}
Consider the many-logic modal structure represented in \cref{A-many-logic-modal-structure-over-L}. Let v be an assignment such that $v_{w_1}(p)=0$, $v_{w_2}(p)=a$ and $v_{w_3}(p)=b$. We have that 

$$v_{w_1}(\Box p)=\bigcap\{\bigcup\{x,0\},\bigcup\{x,0\}\}=x$$

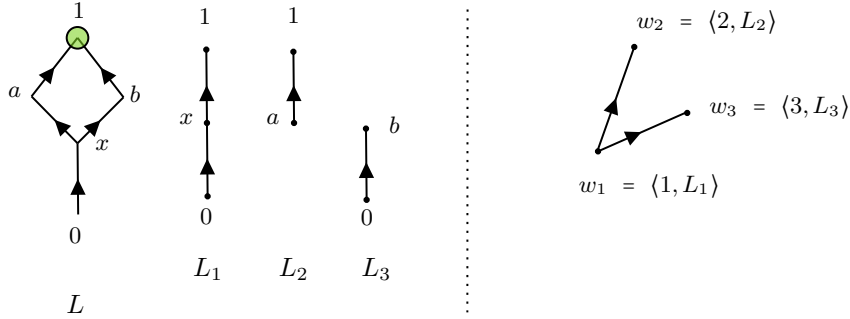
\begin{figure}[H]
\begin{center}
\tikzset{every picture/.style={line width=0.75pt}} 
\begin{tikzpicture}[x=0.75pt,y=0.75pt,yscale=-0.9,xscale=0.9]

    \draw    (124.98,133.6) -- (99.2,107.32) ;
    \draw [shift={(112.09,120.46)}, rotate = 45.56] [fill={rgb, 255:red, 0; green, 0; blue, 0 }  ][line width=0.08]  [draw opacity=0] (8.93,-4.29) -- (0,0) -- (8.93,4.29) -- cycle    ;
    \draw    (124.98,133.6) -- (150.76,107.82) ;
    \draw [shift={(137.87,120.71)}, rotate = 135] [fill={rgb, 255:red, 0; green, 0; blue, 0 }  ][line width=0.08]  [draw opacity=0] (8.93,-4.29) -- (0,0) -- (8.93,4.29) -- cycle    ;
    \draw    (150.76,107.82) -- (124.98,74.47) ;
    \draw [shift={(137.87,91.14)}, rotate = 52.31] [fill={rgb, 255:red, 0; green, 0; blue, 0 }  ][line width=0.08]  [draw opacity=0] (8.93,-4.29) -- (0,0) -- (8.93,4.29) -- cycle    ;
    \draw    (99.2,107.32) -- (124.98,74.47) ;
    \draw [shift={(112.09,90.89)}, rotate = 128.12] [fill={rgb, 255:red, 0; green, 0; blue, 0 }  ][line width=0.08]  [draw opacity=0] (8.93,-4.29) -- (0,0) -- (8.93,4.29) -- cycle    ;
    \draw  [fill={rgb, 255:red, 126; green, 211; blue, 33 }  ,fill opacity=0.58 ] (118.91,74.47) .. controls (118.91,71.12) and (121.63,68.4) .. (124.98,68.4) .. controls (128.33,68.4) and (131.04,71.12) .. (131.04,74.47) .. controls (131.04,77.81) and (128.33,80.53) .. (124.98,80.53) .. controls (121.63,80.53) and (118.91,77.81) .. (118.91,74.47) -- cycle ;
    \draw    (125.34,173.53) -- (124.98,133.6) ;
    \draw [shift={(125.16,153.57)}, rotate = 89.48] [fill={rgb, 255:red, 0; green, 0; blue, 0 }  ][line width=0.08]  [draw opacity=0] (8.93,-4.29) -- (0,0) -- (8.93,4.29) -- cycle    ;
    \draw    (197.87,162.15) -- (197.51,122.21) ;
    \draw [shift={(197.69,142.18)}, rotate = 89.48] [fill={rgb, 255:red, 0; green, 0; blue, 0 }  ][line width=0.08]  [draw opacity=0] (8.93,-4.29) -- (0,0) -- (8.93,4.29) -- cycle    ;
    \draw  [fill={rgb, 255:red, 0; green, 0; blue, 0 }  ,fill opacity=1 ] (196.67,163.35) .. controls (196.67,162.68) and (197.21,162.15) .. (197.87,162.15) .. controls (198.54,162.15) and (199.07,162.68) .. (199.07,163.35) .. controls (199.07,164.01) and (198.54,164.55) .. (197.87,164.55) .. controls (197.21,164.55) and (196.67,164.01) .. (196.67,163.35) -- cycle ;
    \draw  [fill={rgb, 255:red, 0; green, 0; blue, 0 }  ,fill opacity=1 ] (196.3,122.21) .. controls (196.3,121.55) and (196.84,121.01) .. (197.51,121.01) .. controls (198.17,121.01) and (198.71,121.55) .. (198.71,122.21) .. controls (198.71,122.88) and (198.17,123.42) .. (197.51,123.42) .. controls (196.84,123.42) and (196.3,122.88) .. (196.3,122.21) -- cycle ;
    \draw    (197.51,121.01) -- (197.15,81.08) ;
    \draw [shift={(197.33,101.04)}, rotate = 89.48] [fill={rgb, 255:red, 0; green, 0; blue, 0 }  ][line width=0.08]  [draw opacity=0] (8.93,-4.29) -- (0,0) -- (8.93,4.29) -- cycle    ;
    \draw  [fill={rgb, 255:red, 0; green, 0; blue, 0 }  ,fill opacity=1 ] (195.94,81.08) .. controls (195.94,80.41) and (196.48,79.87) .. (197.15,79.87) .. controls (197.81,79.87) and (198.35,80.41) .. (198.35,81.08) .. controls (198.35,81.74) and (197.81,82.28) .. (197.15,82.28) .. controls (196.48,82.28) and (195.94,81.74) .. (195.94,81.08) -- cycle ;
    \draw  [fill={rgb, 255:red, 0; green, 0; blue, 0 }  ,fill opacity=1 ] (245.11,122.21) .. controls (245.11,121.55) and (245.64,121.01) .. (246.31,121.01) .. controls (246.97,121.01) and (247.51,121.55) .. (247.51,122.21) .. controls (247.51,122.88) and (246.97,123.42) .. (246.31,123.42) .. controls (245.64,123.42) and (245.11,122.88) .. (245.11,122.21) -- cycle ;
    \draw    (246.31,121.01) -- (245.95,81.08) ;
    \draw [shift={(246.13,101.04)}, rotate = 89.48] [fill={rgb, 255:red, 0; green, 0; blue, 0 }  ][line width=0.08]  [draw opacity=0] (8.93,-4.29) -- (0,0) -- (8.93,4.29) -- cycle    ;
    \draw  [fill={rgb, 255:red, 0; green, 0; blue, 0 }  ,fill opacity=1 ] (244.74,82.28) .. controls (244.74,81.61) and (245.28,81.08) .. (245.95,81.08) .. controls (246.61,81.08) and (247.15,81.61) .. (247.15,82.28) .. controls (247.15,82.94) and (246.61,83.48) .. (245.95,83.48) .. controls (245.28,83.48) and (244.74,82.94) .. (244.74,82.28) -- cycle ;
    \draw  [fill={rgb, 255:red, 0; green, 0; blue, 0 }  ,fill opacity=1 ] (285.77,165.32) .. controls (285.77,164.65) and (286.31,164.12) .. (286.98,164.12) .. controls (287.64,164.12) and (288.18,164.65) .. (288.18,165.32) .. controls (288.18,165.98) and (287.64,166.52) .. (286.98,166.52) .. controls (286.31,166.52) and (285.77,165.98) .. (285.77,165.32) -- cycle ;
    \draw    (286.98,164.12) -- (286.62,124.18) ;
    \draw [shift={(286.8,144.15)}, rotate = 89.48] [fill={rgb, 255:red, 0; green, 0; blue, 0 }  ][line width=0.08]  [draw opacity=0] (8.93,-4.29) -- (0,0) -- (8.93,4.29) -- cycle    ;
    \draw  [fill={rgb, 255:red, 0; green, 0; blue, 0 }  ,fill opacity=1 ] (285.41,125.39) .. controls (285.41,124.72) and (285.95,124.18) .. (286.62,124.18) .. controls (287.28,124.18) and (287.82,124.72) .. (287.82,125.39) .. controls (287.82,126.05) and (287.28,126.59) .. (286.62,126.59) .. controls (285.95,126.59) and (285.41,126.05) .. (285.41,125.39) -- cycle ;
    \draw    (416.66,136.93) -- (436.97,79.6) ;
    \draw [shift={(426.81,108.26)}, rotate = 109.5] [fill={rgb, 255:red, 0; green, 0; blue, 0 }  ][line width=0.08]  [draw opacity=0] (8.93,-4.29) -- (0,0) -- (8.93,4.29) -- cycle    ;
    \draw  [fill={rgb, 255:red, 0; green, 0; blue, 0 }  ,fill opacity=1 ] (415.46,138.14) .. controls (415.46,137.47) and (416,136.93) .. (416.66,136.93) .. controls (417.33,136.93) and (417.86,137.47) .. (417.86,138.14) .. controls (417.86,138.8) and (417.33,139.34) .. (416.66,139.34) .. controls (416,139.34) and (415.46,138.8) .. (415.46,138.14) -- cycle ;
    \draw  [fill={rgb, 255:red, 0; green, 0; blue, 0 }  ,fill opacity=1 ] (435.76,79.6) .. controls (435.76,78.93) and (436.3,78.39) .. (436.97,78.39) .. controls (437.63,78.39) and (438.17,78.93) .. (438.17,79.6) .. controls (438.17,80.26) and (437.63,80.8) .. (436.97,80.8) .. controls (436.3,80.8) and (435.76,80.26) .. (435.76,79.6) -- cycle ;
    \draw    (416.66,138.14) -- (465.44,116.59) ;
    \draw [shift={(441.05,127.36)}, rotate = 156.17] [fill={rgb, 255:red, 0; green, 0; blue, 0 }  ][line width=0.08]  [draw opacity=0] (8.93,-4.29) -- (0,0) -- (8.93,4.29) -- cycle    ;
    \draw  [fill={rgb, 255:red, 0; green, 0; blue, 0 }  ,fill opacity=1 ] (465.44,116.59) .. controls (465.44,115.92) and (465.97,115.39) .. (466.64,115.39) .. controls (467.3,115.39) and (467.84,115.92) .. (467.84,116.59) .. controls (467.84,117.25) and (467.3,117.79) .. (466.64,117.79) .. controls (465.97,117.79) and (465.44,117.25) .. (465.44,116.59) -- cycle ;
    \draw  [dash pattern={on 0.84pt off 2.51pt}]  (343.55,58.57) -- (343.43,228.44) ;

    \draw (84.36,100.29) node [anchor=north west][inner sep=0.75pt]  [font=\small]  {$a$};
    \draw (152.6,99.28) node [anchor=north west][inner sep=0.75pt]  [font=\small]  {$b$};
    \draw (118.61,179.18) node [anchor=north west][inner sep=0.75pt]  [font=\small]  {$0$};
    \draw (120.68,53.28) node [anchor=north west][inner sep=0.75pt]  [font=\small]  {$1$};
    \draw (191.95,168.6) node [anchor=north west][inner sep=0.75pt]  [font=\small]  {$0$};
    \draw (133.31,129.92) node [anchor=north west][inner sep=0.75pt]  [font=\small]  {$x$};
    \draw (180.48,115.28) node [anchor=north west][inner sep=0.75pt]  [font=\small]  {$x$};
    \draw (191.32,56.54) node [anchor=north west][inner sep=0.75pt]  [font=\small]  {$1$};
    \draw (229.28,115.28) node [anchor=north west][inner sep=0.75pt]  [font=\small]  {$a$};
    \draw (240.93,55.72) node [anchor=north west][inner sep=0.75pt]  [font=\small]  {$1$};
    \draw (282.15,168.97) node [anchor=north west][inner sep=0.75pt]  [font=\small]  {$0$};
    \draw (297.87,116.72) node [anchor=north west][inner sep=0.75pt]  [font=\small]  {$b$};
    \draw (117.04,216.89) node [anchor=north west][inner sep=0.75pt]    {$L$};
    \draw (188.05,196.37) node [anchor=north west][inner sep=0.75pt]    {$L_{1}$};
    \draw (236.04,196.37) node [anchor=north west][inner sep=0.75pt]    {$L_{2}$};
    \draw (282.4,196.37) node [anchor=north west][inner sep=0.75pt]    {$L_{3}$};
    \draw (404.24,149.09) node [anchor=north west][inner sep=0.75pt]  [font=\small]  {$w_{1} \ =\ \langle 1,L_{1} \rangle $};
    \draw (436.47,56.99) node [anchor=north west][inner sep=0.75pt]  [font=\small]  {$w_{2} \ =\ \langle 2,L_{2} \rangle $};
    \draw (477.38,105.16) node [anchor=north west][inner sep=0.75pt]  [font=\small]  {$w_{3} \ =\ \langle 3,L_{3} \rangle $};
\end{tikzpicture}
\end{center}
\caption{A many-logic modal structure over $L$}
\label{A-many-logic-modal-structure-over-L}
\end{figure}
   
\end{ex}

We should point out that many-logic modal structures generalize standard models -- i.e. each Kripke model can be represented by a many-logic modal structures with a common lattice in all worlds. 
For instance, Kripke structures of propositional modal Classical Logic are many-logic modal structures with the two boolean algebra associated to each world.
   
Sometimes it is worth to consider a binary evaluation of formulas in a satisfaction relation. Thus we make use of the filter associated with the base lattice $L$.

\begin{definition}
Let $M = \pair{W, R, s}$ be a many-logic modal $L$-structure and $F$ is the filter associated with $L$. If $w \in W$, we say that $w$ satisfies a formula $\varphi$ (formally, $w \vDash \varphi$) when $v_w(\varphi) \in F$.
\end{definition} 

In next section we will produce some examples where the satisfaction relation contradict NAW in this framework.
   
We finish this section defining a notion of Bisimulation for many-logic modal structures. This is a important notion in traditional modal logic. It is a binary relation between the set of worlds of two Kripke models that identifies two worlds if they have the same behavior (see. \cite[p. 66]{blackburn_modal}). 
Two worlds are bisimilar if the same propositions hold in each world and transition possibilities coincide (i.e. if, in one model, a world $w$ is accessible from the current world, then it must exist an accessible world in the other model bisimilar to $w$). An important result on this subject is that standard modal logic is invariant under bisimulation \cite[p. 67]{blackburn_modal}. The relation is stronger -- for Kripke models of finite image, two worlds are bisimilar if they satisfy the same modal formulas.
    
In our setting one can consider a generalization of bisimulation.

\begin{definition}Let $M$ and $M'$ be two  
many-logic modal structures over a lattice $L$ and let $(W\times W')^=:=\{\pair{\pair{w,A},\pair{w',A'}}\in W\times W' \mid A=A'\}$.
$B\subseteq (W\times W')^=$ is said to be a \emph{bisimulation} if for any $\pair{\pair{w,A},\pair{w',A'}}\in B$ the following conditions hold:

\begin{description}
    \item[Atomic condition] $v_w(p)=v_{w'}(p)$, for any propositional variable $p$;
    \item[Zig] If $\pair{w,A_1} R \pair{v,A_2}$, then there is a world $\pair{v',A_2}$ such that $\pair{w',A_1} R \pair{v',A_2}$ and $\pair{\pair{v,A_2},\pair{v',A_2}}\in B$;
    \item[Zag] If $\pair{w',A_1} R \pair{v',A_2}$, then there is a world $\pair{v,A_2}$ such that $\pair{w,A_1} R \pair{v,A_2}$ and $\pair{\pair{v,A_2},\pair{v',A_2}}\in B$.\footnote{We should point out that, if we have the same sublattice in all worlds, then we obtain the standard notion of bisimulation. Moreover, one may argue that requiring that two bisimilar worlds have the same sublattice is very strong. The study of a weaker condition is a very interesting topic we will leave for future research.}
\end{description} 
\end{definition}

\begin{ex}
Consider the lattice of the previous example (\cref{A-many-logic-modal-structure-over-L}) and the following many-logic modal structures. For simplicity, assume that we have just one variable $p$ and $v_{w}(p)= v_{w'}(p)$ and $v_{v}(p)= v_{v_1'}(p)= v_{v_2'}(p)$. It is easy to see that the relation represented by the red dashed line is a bisimulation.

\begin{figure}[H]
\begin{center}
\tikzset{every picture/.style={line width=0.75pt}} 
\begin{tikzpicture}[x=0.75pt,y=0.75pt,yscale=-1,xscale=1]

    \draw  [fill={rgb, 255:red, 0; green, 0; blue, 0 }  ,fill opacity=1 ] (128.96,131.13) .. controls (128.96,129.91) and (129.95,128.92) .. (131.17,128.92) .. controls (132.39,128.92) and (133.38,129.91) .. (133.38,131.13) .. controls (133.38,132.35) and (132.39,133.34) .. (131.17,133.34) .. controls (129.95,133.34) and (128.96,132.35) .. (128.96,131.13) -- cycle ;
    \draw  [fill={rgb, 255:red, 0; green, 0; blue, 0 }  ,fill opacity=1 ] (232.2,131.13) .. controls (232.2,129.91) and (233.19,128.92) .. (234.42,128.92) .. controls (235.64,128.92) and (236.63,129.91) .. (236.63,131.13) .. controls (236.63,132.35) and (235.64,133.34) .. (234.42,133.34) .. controls (233.19,133.34) and (232.2,132.35) .. (232.2,131.13) -- cycle ;
    \draw    (234.42,131.13) .. controls (256.54,137.03) and (246.22,152.52) .. (235.15,153.26) .. controls (224.64,153.96) and (222.12,146.67) .. (228.22,136.46) ;
    \draw [shift={(229.25,134.82)}, rotate = 123.69] [color={rgb, 255:red, 0; green, 0; blue, 0 }  ][line width=0.75]    (10.93,-3.29) .. controls (6.95,-1.4) and (3.31,-0.3) .. (0,0) .. controls (3.31,0.3) and (6.95,1.4) .. (10.93,3.29)   ;
    \draw    (131.17,131.13) -- (230.2,131.13) ;
    \draw [shift={(232.2,131.13)}, rotate = 180] [color={rgb, 255:red, 0; green, 0; blue, 0 }  ][line width=0.75]    (10.93,-3.29) .. controls (6.95,-1.4) and (3.31,-0.3) .. (0,0) .. controls (3.31,0.3) and (6.95,1.4) .. (10.93,3.29)   ;
    \draw  [fill={rgb, 255:red, 0; green, 0; blue, 0 }  ,fill opacity=1 ] (342.09,130.39) .. controls (342.09,129.17) and (343.08,128.18) .. (344.3,128.18) .. controls (345.52,128.18) and (346.51,129.17) .. (346.51,130.39) .. controls (346.51,131.62) and (345.52,132.61) .. (344.3,132.61) .. controls (343.08,132.61) and (342.09,131.62) .. (342.09,130.39) -- cycle ;
    \draw  [fill={rgb, 255:red, 0; green, 0; blue, 0 }  ,fill opacity=1 ] (471.15,131.13) .. controls (471.15,129.91) and (472.14,128.92) .. (473.36,128.92) .. controls (474.59,128.92) and (475.58,129.91) .. (475.58,131.13) .. controls (475.58,132.35) and (474.59,133.34) .. (473.36,133.34) .. controls (472.14,133.34) and (471.15,132.35) .. (471.15,131.13) -- cycle ;
    \draw    (344.3,130.39) -- (469.15,131.12) ;
    \draw [shift={(471.15,131.13)}, rotate = 180.33] [color={rgb, 255:red, 0; green, 0; blue, 0 }  ][line width=0.75]    (10.93,-3.29) .. controls (6.95,-1.4) and (3.31,-0.3) .. (0,0) .. controls (3.31,0.3) and (6.95,1.4) .. (10.93,3.29)   ;
    \draw    (344.3,130.39) -- (445.12,66.57) ;
    \draw [shift={(446.81,65.5)}, rotate = 147.66] [color={rgb, 255:red, 0; green, 0; blue, 0 }  ][line width=0.75]    (10.93,-3.29) .. controls (6.95,-1.4) and (3.31,-0.3) .. (0,0) .. controls (3.31,0.3) and (6.95,1.4) .. (10.93,3.29)   ;
    \draw  [fill={rgb, 255:red, 0; green, 0; blue, 0 }  ,fill opacity=1 ] (446.08,64.02) .. controls (446.08,62.8) and (447.07,61.81) .. (448.29,61.81) .. controls (449.51,61.81) and (450.5,62.8) .. (450.5,64.02) .. controls (450.5,65.24) and (449.51,66.23) .. (448.29,66.23) .. controls (447.07,66.23) and (446.08,65.24) .. (446.08,64.02) -- cycle ;
    \draw    (448.29,64.02) -- (472.64,127.05) ;
    \draw [shift={(473.36,128.92)}, rotate = 248.88] [color={rgb, 255:red, 0; green, 0; blue, 0 }  ][line width=0.75]    (10.93,-3.29) .. controls (6.95,-1.4) and (3.31,-0.3) .. (0,0) .. controls (3.31,0.3) and (6.95,1.4) .. (10.93,3.29)   ;
    \draw    (473.36,131.13) .. controls (486.37,101.5) and (477.43,71.87) .. (452.07,64.45) ;
    \draw [shift={(450.5,64.02)}, rotate = 14.04] [color={rgb, 255:red, 0; green, 0; blue, 0 }  ][line width=0.75]    (10.93,-3.29) .. controls (6.95,-1.4) and (3.31,-0.3) .. (0,0) .. controls (3.31,0.3) and (6.95,1.4) .. (10.93,3.29)   ;
    \draw [color={rgb, 255:red, 208; green, 2; blue, 27 }  ,draw opacity=1 ] [dash pattern={on 0.84pt off 2.51pt}]  (133.38,131.13) .. controls (158.2,193.19) and (292.17,159.84) .. (339.22,134.84) ;
    \draw [shift={(340.62,134.08)}, rotate = 151.26] [color={rgb, 255:red, 208; green, 2; blue, 27 }  ,draw opacity=1 ][line width=0.75]    (10.93,-3.29) .. controls (6.95,-1.4) and (3.31,-0.3) .. (0,0) .. controls (3.31,0.3) and (6.95,1.4) .. (10.93,3.29)   ;
    \draw [color={rgb, 255:red, 208; green, 2; blue, 27 }  ,draw opacity=1 ] [dash pattern={on 0.84pt off 2.51pt}]  (234.42,131.13) .. controls (278.22,194.65) and (432.62,169.27) .. (472.2,134.4) ;
    \draw [shift={(473.36,133.34)}, rotate = 136.74] [color={rgb, 255:red, 208; green, 2; blue, 27 }  ,draw opacity=1 ][line width=0.75]    (10.93,-3.29) .. controls (6.95,-1.4) and (3.31,-0.3) .. (0,0) .. controls (3.31,0.3) and (6.95,1.4) .. (10.93,3.29)   ;
    \draw [color={rgb, 255:red, 208; green, 2; blue, 27 }  ,draw opacity=1 ] [dash pattern={on 0.84pt off 2.51pt}]  (234.42,131.13) .. controls (227.81,81.97) and (362.85,56.17) .. (444.84,63.9) ;
    \draw [shift={(446.08,64.02)}, rotate = 185.66] [color={rgb, 255:red, 208; green, 2; blue, 27 }  ,draw opacity=1 ][line width=0.75]    (10.93,-3.29) .. controls (6.95,-1.4) and (3.31,-0.3) .. (0,0) .. controls (3.31,0.3) and (6.95,1.4) .. (10.93,3.29)   ;

    \draw (118.82,110.83) node [anchor=north west][inner sep=0.75pt]  [font=\scriptsize]  {$w=\langle 1,L_{1} \rangle $};
    \draw (244.87,110.94) node [anchor=north west][inner sep=0.75pt]  [font=\scriptsize]  {$v=\langle 2,L_{2} \rangle $};
    \draw (341.68,139.41) node [anchor=north west][inner sep=0.75pt]  [font=\scriptsize]  {$w'=\langle 1,L_{1} \rangle $};
    \draw (442.53,40.06) node [anchor=north west][inner sep=0.75pt]  [font=\scriptsize]  {$v_{1} '=\langle 2,L_{2} \rangle $};
    \draw (474.18,139.63) node [anchor=north west][inner sep=0.75pt]  [font=\scriptsize]  {$v_{2} '=\langle 3,L_{2} \rangle $};
    \draw (166.84,188.41) node [anchor=north west][inner sep=0.75pt]    {$M_{1}$};
    \draw (404.31,184.72) node [anchor=north west][inner sep=0.75pt]    {$M_{2}$};
\end{tikzpicture}
\caption{Two bisimilar models}
\end{center}
\end{figure}
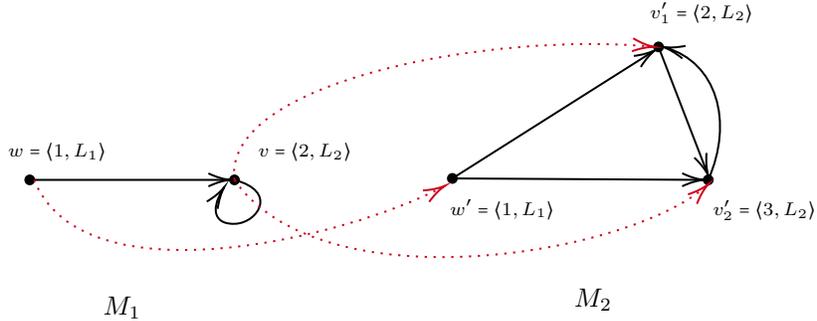

\end{ex}

It is not difficult to show that for any bisimulation $B$, if $\pair{w,w'}\in B$, then for any modal formula $\varphi$, $v_w(\varphi)=v_{w'}(\varphi)$. The base case holds for the atomic condition. It is also easy to verify that the induction steps for the Boolean connectives  hold. The hard part is the induction step for formulas of type $\Box \varphi$. From the Zig condition, if $wRu$ there is a $u'$ such that $w'R'u'$ and $uBu'$. Hence, by induction hypothesis, $\displaystyle \big\{\bigcap\{ x\in L_w \mid x\leq v_u(\varphi)\} \mid w R u \big\} \subseteq \big\{\bigcap\{ x\in L_{w'} \mid x\leq v_{u'}(\varphi))\} \mid w' R u' \big\} $. Hence, $v_w(\Box \varphi)\leq v_{w'}(\Box \varphi)$. The Zag condition, similarly ensures $v_w(\Box \varphi)\geq v_{w'}(\Box \varphi)$. And so, $v_w(\Box \varphi)= v_{w'}(\Box \varphi)$.

\section{Accessibility and necessity}

With the framework of many-logic modal structure, we can now study some natural properties and possibilities. We start discussing the locality of axiom $K$ (we assume material implication is defined as $\varphi \to \psi := \lnot \varphi \lor \psi$).

Since in our approach worlds can have different lattices, it is easy to find examples where axiom $K$ fails to hold. Moreover, even if the local lattice is well behaved, $K$ can fail. 

\begin{ex}  Let $w = \pair{1, L_1}$ and $w' = \pair{2, L_2}$ be worlds built with lattices $L_1$ and $L_2$ in \cref{figure-fail-k} (complement function is represented with red arrows). Let $R = \{\pair{w,w'}\}$, $W = \{w, w'\}$ and $M = \pair{W, R, v}$ such that $v_w(p) = a $, $v_w(q) = b $ and $v_{w'}(p) = a_1$ $v_{w'}(q) = b$ for propositional variables $p$ and $q$.

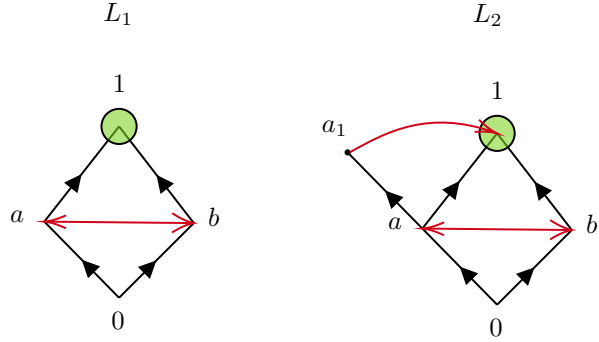
\begin{figure}[H]
\begin{center}
\tikzset{every picture/.style={line width=0.75pt}} 
\begin{tikzpicture}[x=0.75pt,y=0.75pt,yscale=-1,xscale=1]

    \draw    (201.29,185.93) -- (163.59,147.49) ;
    \draw [shift={(182.44,166.71)}, rotate = 45.56] [fill={rgb, 255:red, 0; green, 0; blue, 0 }  ][line width=0.08]  [draw opacity=0] (8.93,-4.29) -- (0,0) -- (8.93,4.29) -- cycle    ;
    \draw    (201.29,185.93) -- (238.98,148.23) ;
    \draw [shift={(220.13,167.08)}, rotate = 135] [fill={rgb, 255:red, 0; green, 0; blue, 0 }  ][line width=0.08]  [draw opacity=0] (8.93,-4.29) -- (0,0) -- (8.93,4.29) -- cycle    ;
    \draw    (238.98,148.23) -- (201.29,99.45) ;
    \draw [shift={(220.13,123.84)}, rotate = 52.31] [fill={rgb, 255:red, 0; green, 0; blue, 0 }  ][line width=0.08]  [draw opacity=0] (8.93,-4.29) -- (0,0) -- (8.93,4.29) -- cycle    ;
    \draw    (163.59,147.49) -- (201.29,99.45) ;
    \draw [shift={(182.44,123.47)}, rotate = 128.12] [fill={rgb, 255:red, 0; green, 0; blue, 0 }  ][line width=0.08]  [draw opacity=0] (8.93,-4.29) -- (0,0) -- (8.93,4.29) -- cycle    ;
    \draw  [fill={rgb, 255:red, 126; green, 211; blue, 33 }  ,fill opacity=0.58 ] (192.42,99.45) .. controls (192.42,94.56) and (196.39,90.58) .. (201.29,90.58) .. controls (206.18,90.58) and (210.15,94.56) .. (210.15,99.45) .. controls (210.15,104.35) and (206.18,108.32) .. (201.29,108.32) .. controls (196.39,108.32) and (192.42,104.35) .. (192.42,99.45) -- cycle ;
    \draw    (391.95,189.43) -- (354.26,150.99) ;
    \draw [shift={(373.1,170.21)}, rotate = 45.56] [fill={rgb, 255:red, 0; green, 0; blue, 0 }  ][line width=0.08]  [draw opacity=0] (8.93,-4.29) -- (0,0) -- (8.93,4.29) -- cycle    ;
    \draw    (391.95,189.43) -- (429.64,151.73) ;
    \draw [shift={(410.8,170.58)}, rotate = 135] [fill={rgb, 255:red, 0; green, 0; blue, 0 }  ][line width=0.08]  [draw opacity=0] (8.93,-4.29) -- (0,0) -- (8.93,4.29) -- cycle    ;
    \draw    (429.64,151.73) -- (391.95,102.95) ;
    \draw [shift={(410.8,127.34)}, rotate = 52.31] [fill={rgb, 255:red, 0; green, 0; blue, 0 }  ][line width=0.08]  [draw opacity=0] (8.93,-4.29) -- (0,0) -- (8.93,4.29) -- cycle    ;
    \draw    (354.26,150.99) -- (391.95,102.95) ;
    \draw [shift={(373.1,126.97)}, rotate = 128.12] [fill={rgb, 255:red, 0; green, 0; blue, 0 }  ][line width=0.08]  [draw opacity=0] (8.93,-4.29) -- (0,0) -- (8.93,4.29) -- cycle    ;
    \draw  [fill={rgb, 255:red, 126; green, 211; blue, 33 }  ,fill opacity=0.58 ] (383.08,102.95) .. controls (383.08,98.05) and (387.05,94.08) .. (391.95,94.08) .. controls (396.85,94.08) and (400.82,98.05) .. (400.82,102.95) .. controls (400.82,107.85) and (396.85,111.82) .. (391.95,111.82) .. controls (387.05,111.82) and (383.08,107.85) .. (383.08,102.95) -- cycle ;
    \draw    (354.26,150.99) -- (316.56,112.56) ;
    \draw [shift={(335.41,131.78)}, rotate = 45.56] [fill={rgb, 255:red, 0; green, 0; blue, 0 }  ][line width=0.08]  [draw opacity=0] (8.93,-4.29) -- (0,0) -- (8.93,4.29) -- cycle    ;
    \draw  [fill={rgb, 255:red, 0; green, 0; blue, 0 }  ,fill opacity=1 ] (315.51,112.56) .. controls (315.51,111.98) and (315.98,111.51) .. (316.56,111.51) .. controls (317.14,111.51) and (317.61,111.98) .. (317.61,112.56) .. controls (317.61,113.14) and (317.14,113.61) .. (316.56,113.61) .. controls (315.98,113.61) and (315.51,113.14) .. (315.51,112.56) -- cycle ;
    \draw [color={rgb, 255:red, 208; green, 2; blue, 27 }  ,draw opacity=1 ]   (165.59,147.51) -- (236.98,148.21) ;
    \draw [shift={(238.98,148.23)}, rotate = 180.56] [color={rgb, 255:red, 208; green, 2; blue, 27 }  ,draw opacity=1 ][line width=0.75]    (10.93,-4.9) .. controls (6.95,-2.3) and (3.31,-0.67) .. (0,0) .. controls (3.31,0.67) and (6.95,2.3) .. (10.93,4.9)   ;
    \draw [shift={(163.59,147.49)}, rotate = 0.56] [color={rgb, 255:red, 208; green, 2; blue, 27 }  ,draw opacity=1 ][line width=0.75]    (10.93,-3.29) .. controls (6.95,-1.4) and (3.31,-0.3) .. (0,0) .. controls (3.31,0.3) and (6.95,1.4) .. (10.93,3.29)   ;
    \draw [color={rgb, 255:red, 208; green, 2; blue, 27 }  ,draw opacity=1 ]   (356.26,151.01) -- (427.64,151.71) ;
    \draw [shift={(429.64,151.73)}, rotate = 180.56] [color={rgb, 255:red, 208; green, 2; blue, 27 }  ,draw opacity=1 ][line width=0.75]    (10.93,-4.9) .. controls (6.95,-2.3) and (3.31,-0.67) .. (0,0) .. controls (3.31,0.67) and (6.95,2.3) .. (10.93,4.9)   ;
    \draw [shift={(354.26,150.99)}, rotate = 0.56] [color={rgb, 255:red, 208; green, 2; blue, 27 }  ,draw opacity=1 ][line width=0.75]    (10.93,-3.29) .. controls (6.95,-1.4) and (3.31,-0.3) .. (0,0) .. controls (3.31,0.3) and (6.95,1.4) .. (10.93,3.29)   ;
    \draw [color={rgb, 255:red, 208; green, 2; blue, 27 }  ,draw opacity=1 ]   (317.61,112.56) .. controls (336.23,102) and (357.73,90.82) .. (390.45,102.41) ;
    \draw [shift={(391.95,102.95)}, rotate = 200.39] [color={rgb, 255:red, 208; green, 2; blue, 27 }  ,draw opacity=1 ][line width=0.75]    (10.93,-3.29) .. controls (6.95,-1.4) and (3.31,-0.3) .. (0,0) .. controls (3.31,0.3) and (6.95,1.4) .. (10.93,3.29)   ;

    \draw (144.9,140.57) node [anchor=north west][inner sep=0.75pt]    {$a$};
    \draw (244.67,139.09) node [anchor=north west][inner sep=0.75pt]    {$b$};
    \draw (195.89,191.57) node [anchor=north west][inner sep=0.75pt]    {$0$};
    \draw (196.63,71.84) node [anchor=north west][inner sep=0.75pt]    {$1$};
    \draw (335.56,144.07) node [anchor=north west][inner sep=0.75pt]    {$a$};
    \draw (435.34,142.59) node [anchor=north west][inner sep=0.75pt]    {$b$};
    \draw (386.56,195.07) node [anchor=north west][inner sep=0.75pt]    {$0$};
    \draw (387.3,75.33) node [anchor=north west][inner sep=0.75pt]    {$1$};
    \draw (191.89,35.46) node [anchor=north west][inner sep=0.75pt]    {$L_{1}$};
    \draw (378.18,36.33) node [anchor=north west][inner sep=0.75pt]    {$L_{2}$};
    \draw (302.07,94.97) node [anchor=north west][inner sep=0.75pt]    {$a_{1}$};
\end{tikzpicture}
\end{center}
\caption{Odd negation for a value}
\label{figure-fail-k}
\end{figure}

It is easy to see that $v_w(\Box (p\to q))=1$ and $v_w(\Box p\to \Box q)= b$. So, $K$ does not hold at $w$. Note that, in this case, the negation has a rigid interpretation. 

\end{ex}

However, in relatively simple/weak global conditions, the verification of the axiom $K$ only depends on local conditions (world by world). In the evaluation of $\Box (p\to q)\to (\Box p\to \Box q)$ in a given world, only $p\to q$ is calculated in the accessible worlds, all other calculations are made in the world where we are evaluating $K$, using the value of $p$ and $q$ in the accessible worlds only as parameters.
Thus, as we have defined $a \to b:=-a+b$, if the  operation ``$-$" is well behaved (as we just see, being rigid is not enough), the verification of $K$ only depends on what happens in the current world, namely properties of the local lattice. Analogously, we can argue that necessity rule also only depends on local properties.\footnote{The  full characterization of the validity of axiom $K$ and rule of necessitation will be the studied in a more technical article in preparation.}

Next we will point out some other interesting facts regarding evaluation of necessity in many-logic modal structures. For example,  there are natural situations in which a many-logic modal structure is such that there are worlds for which (i) a formula is not necessary even though it is considered true in all accessed worlds and (ii) a formula is necessary even though it is considered false in some of the accessed worlds.
As this phenomenon confronts what usually happens in modal logic, it may receive some skeptical looks. However, bearing in mind that our approach is based on the assumption that the concept of truth can vary from world to world, it is natural that, when we relativize what happens in worlds accessible to the current world, unexpected conclusions may follow depending on how we agree on the way worlds in different logics should transfer information. 
\begin{ex}
Let us consider the trivial lattice $L$ with a sublattice $A$:

\begin{figure}[H]
\begin{center}
\tikzset{every picture/.style={line width=0.75pt}} 
\begin{tikzpicture}[x=0.75pt,y=0.75pt,yscale=-1,xscale=1]

	\draw    (150.96,73.62) -- (150.96,25.02) ;
	\draw [shift={(150.96,49.32)}, rotate = 90] [fill={rgb, 255:red, 0; green, 0; blue, 0 }  ][line width=0.08]  [draw opacity=0] (8.93,-4.29) -- (0,0) -- (8.93,4.29) -- cycle    ;
	\draw    (150.56,111) -- (150.96,73.62) ;
	\draw [shift={(150.76,92.31)}, rotate = 90.62] [fill={rgb, 255:red, 0; green, 0; blue, 0 }  ][line width=0.08]  [draw opacity=0] (8.93,-4.29) -- (0,0) -- (8.93,4.29) -- cycle    ;
	\draw  [fill={rgb, 255:red, 126; green, 211; blue, 33 }  ,fill opacity=0.58 ] (146.12,76.63) .. controls (146.12,73.88) and (148.29,71.65) .. (150.96,71.65) .. controls (153.64,71.65) and (155.81,73.88) .. (155.81,76.63) .. controls (155.81,79.39) and (153.64,81.62) .. (150.96,81.62) .. controls (148.29,81.62) and (146.12,79.39) .. (146.12,76.63) -- cycle ;
	\draw  [fill={rgb, 255:red, 126; green, 211; blue, 33 }  ,fill opacity=0.58 ] (146.12,25.02) .. controls (146.12,22.27) and (148.29,20.04) .. (150.96,20.04) .. controls (153.64,20.04) and (155.81,22.27) .. (155.81,25.02) .. controls (155.81,27.77) and (153.64,30) .. (150.96,30) .. controls (148.29,30) and (146.12,27.77) .. (146.12,25.02) -- cycle ;
	\draw    (257.18,112.45) -- (257.59,26.27) ;
	\draw [shift={(257.39,69.36)}, rotate = 90.27] [fill={rgb, 255:red, 0; green, 0; blue, 0 }  ][line width=0.08]  [draw opacity=0] (8.93,-4.29) -- (0,0) -- (8.93,4.29) -- cycle    ;
	\draw  [fill={rgb, 255:red, 126; green, 211; blue, 33 }  ,fill opacity=0.58 ] (252.74,26.27) .. controls (252.74,23.51) and (254.91,21.28) .. (257.59,21.28) .. controls (260.27,21.28) and (262.44,23.51) .. (262.44,26.27) .. controls (262.44,29.02) and (260.27,31.25) .. (257.59,31.25) .. controls (254.91,31.25) and (252.74,29.02) .. (252.74,26.27) -- cycle ;

	\draw (145.96,3.79) node [anchor=north west][inner sep=0.75pt]  [font=\small]  {$1$};
	\draw (160.1,61.6) node [anchor=north west][inner sep=0.75pt]  [font=\small]  {$1'$};
	\draw (145.16,116.8) node [anchor=north west][inner sep=0.75pt]  [font=\small]  {$0$};
	\draw (156.46,111.88) node [anchor=north west][inner sep=0.75pt]  [font=\scriptsize]  {$( -1,\ -1')$};
	\draw (121.53,11.61) node [anchor=north west][inner sep=0.75pt]  [font=\normalsize]  {$L$};
	\draw (329.73,18.43) node [anchor=north west][inner sep=0.75pt]  [font=\small]  {$Filter\ =\{1,1'\}$};
	\draw (218.95,13.92) node [anchor=north west][inner sep=0.75pt]  [font=\normalsize]  {$A$};
	\draw (252.59,5.04) node [anchor=north west][inner sep=0.75pt]  [font=\small]  {$1$};
	\draw (251.78,115.05) node [anchor=north west][inner sep=0.75pt]  [font=\small]  {$0$};
	\draw (268.3,110.91) node [anchor=north west][inner sep=0.75pt]  [font=\scriptsize]  {$( -1)$};
\end{tikzpicture}
\end{center}
\caption{Three valued lattice with two values in the filter.}
\end{figure}
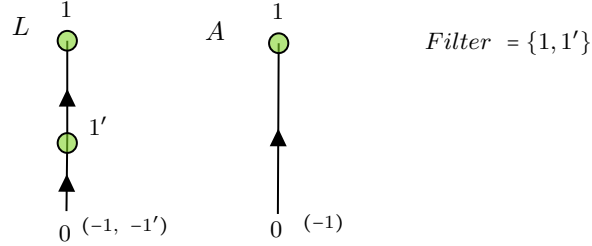
We will consider the situation where we have the worlds $w = \pair{1, A}$ and $w' = \pair{2, L}$. Let $R = \{\pair{w,w'}, \pair{w,w}\}$, $W = \{w, w'\}$ and $M = \pair{W, R, v}$ such that $v_w(\alpha) = 1$ and $v_{w'}(\alpha) = 1'$ for a propositional variable $\alpha$. 

Now we observe that $w \vDash \alpha$, $w \nvDash \lnot \alpha$, $w' \vDash \alpha$, $w' \nvDash \lnot \alpha$. So both worlds accessed by $w$ validate $\alpha$ and both do not validate $\lnot \alpha$. Nonetheless, $w \nvDash \Box \alpha$, for $v_w(\Box \alpha) = \bigcap_A \{1, 1'\} = 1 \cap 0 = 0$.
\end{ex}

This lattice thus produces an interesting example, where a formula can be valid in all accessed worlds while the formula is not necessary. In the example, the phenomenon results from the fact that, once one considers the value $1'$ from the point of view of $w$, the value can only be understood as the max value (in $A$) bellow $1'$. The formula $\alpha$ is indeed true for $w'$, but its truth-value is not directly understandable for $w$. And it is this opening that renders this interesting case possible.

Let us now consider a slightly more complicated case, where a formula can be false in an accessed world and still be necessary. 

\begin{ex}
Consider the following lattice $L$ with sublaticces $A_1$ and $A_2$:

\begin{figure}[H]
\begin{center}
\tikzset{every picture/.style={line width=0.75pt}} 
\begin{tikzpicture}[x=0.75pt,y=0.75pt,yscale=-1,xscale=1]

	\draw    (187.09,190.25) -- (217.02,165.11) ;
	\draw [shift={(202.05,177.68)}, rotate = 139.97] [fill={rgb, 255:red, 0; green, 0; blue, 0 }  ][line width=0.08]  [draw opacity=0] (8.93,-4.29) -- (0,0) -- (8.93,4.29) -- cycle    ;
	\draw    (187.09,190.25) -- (157.16,165.11) ;
	\draw [shift={(172.13,177.68)}, rotate = 40.03] [fill={rgb, 255:red, 0; green, 0; blue, 0 }  ][line width=0.08]  [draw opacity=0] (8.93,-4.29) -- (0,0) -- (8.93,4.29) -- cycle    ;
	\draw    (217.02,165.11) -- (187.69,136.38) ;
	\draw [shift={(202.35,150.75)}, rotate = 44.41] [fill={rgb, 255:red, 0; green, 0; blue, 0 }  ][line width=0.08]  [draw opacity=0] (8.93,-4.29) -- (0,0) -- (8.93,4.29) -- cycle    ;
	\draw    (157.16,165.11) -- (187.69,136.38) ;
	\draw [shift={(172.43,150.75)}, rotate = 136.74] [fill={rgb, 255:red, 0; green, 0; blue, 0 }  ][line width=0.08]  [draw opacity=0] (8.93,-4.29) -- (0,0) -- (8.93,4.29) -- cycle    ;
	\draw    (187.69,136.38) -- (157.16,105.26) ;
	\draw [shift={(172.43,120.82)}, rotate = 45.56] [fill={rgb, 255:red, 0; green, 0; blue, 0 }  ][line width=0.08]  [draw opacity=0] (8.93,-4.29) -- (0,0) -- (8.93,4.29) -- cycle    ;
	\draw    (187.69,136.38) -- (218.22,105.86) ;
	\draw [shift={(202.95,121.12)}, rotate = 135] [fill={rgb, 255:red, 0; green, 0; blue, 0 }  ][line width=0.08]  [draw opacity=0] (8.93,-4.29) -- (0,0) -- (8.93,4.29) -- cycle    ;
	\draw    (218.22,105.86) -- (187.69,66.35) ;
	\draw [shift={(202.95,86.11)}, rotate = 52.31] [fill={rgb, 255:red, 0; green, 0; blue, 0 }  ][line width=0.08]  [draw opacity=0] (8.93,-4.29) -- (0,0) -- (8.93,4.29) -- cycle    ;
	\draw    (157.16,105.26) -- (187.69,66.35) ;
	\draw [shift={(172.43,85.81)}, rotate = 128.12] [fill={rgb, 255:red, 0; green, 0; blue, 0 }  ][line width=0.08]  [draw opacity=0] (8.93,-4.29) -- (0,0) -- (8.93,4.29) -- cycle    ;
	\draw  [fill={rgb, 255:red, 126; green, 211; blue, 33 }  ,fill opacity=0.58 ] (180.51,66.35) .. controls (180.51,62.39) and (183.72,59.17) .. (187.69,59.17) .. controls (191.66,59.17) and (194.87,62.39) .. (194.87,66.35) .. controls (194.87,70.32) and (191.66,73.54) .. (187.69,73.54) .. controls (183.72,73.54) and (180.51,70.32) .. (180.51,66.35) -- cycle ;
	\draw  [fill={rgb, 255:red, 126; green, 211; blue, 33 }  ,fill opacity=0.58 ] (149.98,105.26) .. controls (149.98,101.29) and (153.2,98.08) .. (157.16,98.08) .. controls (161.13,98.08) and (164.35,101.29) .. (164.35,105.26) .. controls (164.35,109.23) and (161.13,112.44) .. (157.16,112.44) .. controls (153.2,112.44) and (149.98,109.23) .. (149.98,105.26) -- cycle ;
	\draw    (368.94,183.84) -- (396.92,160.34) ;
	\draw [shift={(382.93,172.09)}, rotate = 139.97] [fill={rgb, 255:red, 0; green, 0; blue, 0 }  ][line width=0.08]  [draw opacity=0] (8.93,-4.29) -- (0,0) -- (8.93,4.29) -- cycle    ;
	\draw    (368.94,183.84) -- (340.96,160.34) ;
	\draw [shift={(354.95,172.09)}, rotate = 40.03] [fill={rgb, 255:red, 0; green, 0; blue, 0 }  ][line width=0.08]  [draw opacity=0] (8.93,-4.29) -- (0,0) -- (8.93,4.29) -- cycle    ;
	\draw    (396.92,160.34) -- (340.96,104.39) ;
	\draw [shift={(368.94,132.37)}, rotate = 45] [fill={rgb, 255:red, 0; green, 0; blue, 0 }  ][line width=0.08]  [draw opacity=0] (8.93,-4.29) -- (0,0) -- (8.93,4.29) -- cycle    ;
	\draw    (340.96,160.34) -- (340.96,104.39) ;
	\draw [shift={(340.96,132.37)}, rotate = 90] [fill={rgb, 255:red, 0; green, 0; blue, 0 }  ][line width=0.08]  [draw opacity=0] (8.93,-4.29) -- (0,0) -- (8.93,4.29) -- cycle    ;
	\draw    (340.96,104.39) -- (369.5,68.02) ;
	\draw [shift={(355.23,86.2)}, rotate = 128.12] [fill={rgb, 255:red, 0; green, 0; blue, 0 }  ][line width=0.08]  [draw opacity=0] (8.93,-4.29) -- (0,0) -- (8.93,4.29) -- cycle    ;
	\draw  [fill={rgb, 255:red, 126; green, 211; blue, 33 }  ,fill opacity=0.58 ] (362.78,68.02) .. controls (362.78,64.31) and (365.79,61.31) .. (369.5,61.31) .. controls (373.21,61.31) and (376.21,64.31) .. (376.21,68.02) .. controls (376.21,71.73) and (373.21,74.73) .. (369.5,74.73) .. controls (365.79,74.73) and (362.78,71.73) .. (362.78,68.02) -- cycle ;
	\draw  [fill={rgb, 255:red, 126; green, 211; blue, 33 }  ,fill opacity=0.58 ] (334.25,104.39) .. controls (334.25,100.68) and (337.25,97.67) .. (340.96,97.67) .. controls (344.67,97.67) and (347.68,100.68) .. (347.68,104.39) .. controls (347.68,108.1) and (344.67,111.1) .. (340.96,111.1) .. controls (337.25,111.1) and (334.25,108.1) .. (334.25,104.39) -- cycle ;
	\draw    (489.54,184.96) -- (517.52,161.46) ;
	\draw [shift={(503.53,173.21)}, rotate = 139.97] [fill={rgb, 255:red, 0; green, 0; blue, 0 }  ][line width=0.08]  [draw opacity=0] (8.93,-4.29) -- (0,0) -- (8.93,4.29) -- cycle    ;
	\draw    (489.54,184.96) -- (461.57,161.46) ;
	\draw [shift={(475.56,173.21)}, rotate = 40.03] [fill={rgb, 255:red, 0; green, 0; blue, 0 }  ][line width=0.08]  [draw opacity=0] (8.93,-4.29) -- (0,0) -- (8.93,4.29) -- cycle    ;
	\draw    (517.52,161.46) -- (518.64,106.07) ;
	\draw [shift={(518.08,133.76)}, rotate = 91.16] [fill={rgb, 255:red, 0; green, 0; blue, 0 }  ][line width=0.08]  [draw opacity=0] (8.93,-4.29) -- (0,0) -- (8.93,4.29) -- cycle    ;
	\draw    (461.57,161.46) -- (518.64,106.07) ;
	\draw [shift={(490.1,133.76)}, rotate = 135.86] [fill={rgb, 255:red, 0; green, 0; blue, 0 }  ][line width=0.08]  [draw opacity=0] (8.93,-4.29) -- (0,0) -- (8.93,4.29) -- cycle    ;
	\draw    (518.64,106.07) -- (490.1,69.14) ;
	\draw [shift={(504.37,87.6)}, rotate = 52.31] [fill={rgb, 255:red, 0; green, 0; blue, 0 }  ][line width=0.08]  [draw opacity=0] (8.93,-4.29) -- (0,0) -- (8.93,4.29) -- cycle    ;
	\draw  [fill={rgb, 255:red, 126; green, 211; blue, 33 }  ,fill opacity=0.58 ] (483.39,69.14) .. controls (483.39,65.43) and (486.4,62.42) .. (490.1,62.42) .. controls (493.81,62.42) and (496.82,65.43) .. (496.82,69.14) .. controls (496.82,72.85) and (493.81,75.85) .. (490.1,75.85) .. controls (486.4,75.85) and (483.39,72.85) .. (483.39,69.14) -- cycle ;
	\draw  [dash pattern={on 0.84pt off 2.51pt}]  (283.62,64.38) -- (283.62,181.66) ;

	\draw (181.76,39.8) node [anchor=north west][inner sep=0.75pt]    {$1$};
	\draw (181.56,191.1) node [anchor=north west][inner sep=0.75pt]    {$0$};
	\draw (140.46,157.58) node [anchor=north west][inner sep=0.75pt]    {$e$};
	\draw (120.97,57.22) node [anchor=north west][inner sep=0.75pt]  [font=\normalsize]  {$L$};
	\draw (202.83,43.32) node [anchor=north west][inner sep=0.75pt]  [font=\footnotesize]  {$F\ =\{1,a\}$};
	\draw (135.67,96.53) node [anchor=north west][inner sep=0.75pt]    {$a$};
	\draw (221.86,95.93) node [anchor=north west][inner sep=0.75pt]    {$b$};
	\draw (201.12,129.72) node [anchor=north west][inner sep=0.75pt]  [font=\footnotesize]  {$c$};
	\draw (220.66,158.78) node [anchor=north west][inner sep=0.75pt]    {$f$};
	\draw (363.57,42.55) node [anchor=north west][inner sep=0.75pt]    {$1$};
	\draw (363.45,184.14) node [anchor=north west][inner sep=0.75pt]    {$0$};
	\draw (324.96,152.81) node [anchor=north west][inner sep=0.75pt]    {$e$};
	\draw (312.61,53.95) node [anchor=north west][inner sep=0.75pt]  [font=\normalsize]  {$A_{1}$};
	\draw (320.48,95.74) node [anchor=north west][inner sep=0.75pt]    {$a$};
	\draw (399.93,153.93) node [anchor=north west][inner sep=0.75pt]    {$f$};
	\draw (484.17,43.67) node [anchor=north west][inner sep=0.75pt]    {$1$};
	\draw (484.05,185.26) node [anchor=north west][inner sep=0.75pt]    {$0$};
	\draw (445.56,153.93) node [anchor=north west][inner sep=0.75pt]    {$e$};
	\draw (521.66,96.3) node [anchor=north west][inner sep=0.75pt]    {$b$};
	\draw (520.54,155.05) node [anchor=north west][inner sep=0.75pt]    {$f$};
	\draw (434.24,55.44) node [anchor=north west][inner sep=0.75pt]  [font=\small]  {$A_{2}$};
\end{tikzpicture}
\end{center}
\caption{Lattices with necessary though false in accessed world.}
\label{lattice-contain-up-down-lattices}
\end{figure}
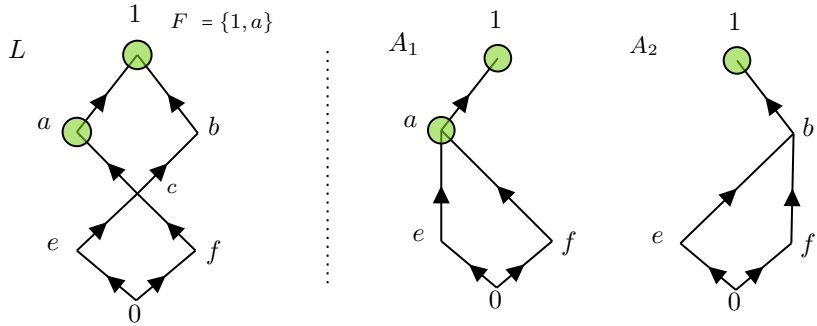
Now we build the following model:
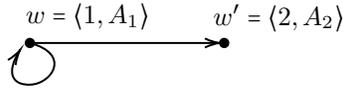
\begin{figure}[H]
\begin{center}
    \tikzset{every picture/.style={line width=0.75pt}} 
    \begin{tikzpicture}[x=0.75pt,y=0.75pt,yscale=-.7,xscale=.7]
        
        \draw  [fill={rgb, 255:red, 0; green, 0; blue, 0 }  ,fill opacity=1 ] (198,140.5) .. controls (198,138.84) and (199.34,137.5) .. (201,137.5) .. controls (202.66,137.5) and (204,138.84) .. (204,140.5) .. controls (204,142.16) and (202.66,143.5) .. (201,143.5) .. controls (199.34,143.5) and (198,142.16) .. (198,140.5) -- cycle ;
        \draw  [fill={rgb, 255:red, 0; green, 0; blue, 0 }  ,fill opacity=1 ] (338,140.5) .. controls (338,138.84) and (339.34,137.5) .. (341,137.5) .. controls (342.66,137.5) and (344,138.84) .. (344,140.5) .. controls (344,142.16) and (342.66,143.5) .. (341,143.5) .. controls (339.34,143.5) and (338,142.16) .. (338,140.5) -- cycle ;
        \draw    (201,140.5) .. controls (231,148.5) and (217,169.5) .. (202,170.5) .. controls (187.52,171.47) and (184.22,161.26) .. (193,147.06) ;
        \draw [shift={(194,145.5)}, rotate = 483.69] [color={rgb, 255:red, 0; green, 0; blue, 0 }  ][line width=0.75]    (10.93,-3.29) .. controls (6.95,-1.4) and (3.31,-0.3) .. (0,0) .. controls (3.31,0.3) and (6.95,1.4) .. (10.93,3.29)   ;
        \draw    (201,140.5) -- (336,140.5) ;
        \draw [shift={(338,140.5)}, rotate = 180] [color={rgb, 255:red, 0; green, 0; blue, 0 }  ][line width=0.75]    (10.93,-3.29) .. controls (6.95,-1.4) and (3.31,-0.3) .. (0,0) .. controls (3.31,0.3) and (6.95,1.4) .. (10.93,3.29)   ;
        
        \draw (196,109.9) node [anchor=north west][inner sep=0.75pt]    {$w=\langle 1,A_1\rangle $};
        \draw (331,110.9) node [anchor=north west][inner sep=0.75pt]    {$w'=\langle 2,A_2\rangle $};
    \end{tikzpicture}
\end{center}
\caption{Model with necessary though false in accessible}
\end{figure}

We attribute $v_w(\alpha) = 1$ and $v_{w'}(\alpha) = b$, for a propositional variable $\alpha$. In this case we have $w \vDash \alpha$ and $w' \nvDash \alpha$. Moreover, $w \vDash \Box \alpha$, since  $v_w(\Box \alpha)=a$.
\end{ex}

Note, however, that $\Box \alpha$ is true in $w$ even though it accesses a world in which $\alpha$ fails. This comes from the fact that $w'$'s value for $\alpha$ is `perceived' as `big enough' in $w$'s universe.

Next we present examples based on lattices related to logics widely studied in the literature.

\begin{ex} Consider the example with simple Heyting structure and double negation.

\begin{figure}[H]
\begin{center}
\tikzset{every picture/.style={line width=0.75pt}} 
\begin{tikzpicture}[x=0.75pt,y=0.75pt,yscale=-1,xscale=1]

	\draw    (244.16,151.15) -- (243.36,109.84) ;
	\draw [shift={(243.76,130.49)}, rotate = 88.89] [fill={rgb, 255:red, 0; green, 0; blue, 0 }  ][line width=0.08]  [draw opacity=0] (8.93,-4.29) -- (0,0) -- (8.93,4.29) -- cycle    ;
	\draw    (243.36,109.84) -- (242.55,58.89) ;
	\draw [shift={(242.95,84.36)}, rotate = 89.1] [fill={rgb, 255:red, 0; green, 0; blue, 0 }  ][line width=0.08]  [draw opacity=0] (8.93,-4.29) -- (0,0) -- (8.93,4.29) -- cycle    ;
	\draw  [fill={rgb, 255:red, 126; green, 211; blue, 33 }  ,fill opacity=0.58 ] (232.93,58.89) .. controls (232.93,53.58) and (237.24,49.27) .. (242.55,49.27) .. controls (247.87,49.27) and (252.18,53.58) .. (252.18,58.89) .. controls (252.18,64.21) and (247.87,68.52) .. (242.55,68.52) .. controls (237.24,68.52) and (232.93,64.21) .. (232.93,58.89) -- cycle ;
	\draw    (346.85,150.35) -- (346.85,66.11) ;
	\draw [shift={(346.85,108.23)}, rotate = 90] [fill={rgb, 255:red, 0; green, 0; blue, 0 }  ][line width=0.08]  [draw opacity=0] (8.93,-4.29) -- (0,0) -- (8.93,4.29) -- cycle    ;
	\draw  [fill={rgb, 255:red, 126; green, 211; blue, 33 }  ,fill opacity=0.58 ] (337.22,59.69) .. controls (337.22,54.38) and (341.53,50.07) .. (346.85,50.07) .. controls (352.16,50.07) and (356.47,54.38) .. (356.47,59.69) .. controls (356.47,65.01) and (352.16,69.32) .. (346.85,69.32) .. controls (341.53,69.32) and (337.22,65.01) .. (337.22,59.69) -- cycle ;
	\draw  [fill={rgb, 255:red, 0; green, 0; blue, 0 }  ,fill opacity=1 ] (240.55,109.84) .. controls (240.55,108.29) and (241.81,107.03) .. (243.36,107.03) .. controls (244.91,107.03) and (246.16,108.29) .. (246.16,109.84) .. controls (246.16,111.39) and (244.91,112.64) .. (243.36,112.64) .. controls (241.81,112.64) and (240.55,111.39) .. (240.55,109.84) -- cycle ;

	\draw (279,-14.6) node [anchor=north west][inner sep=0.75pt]    {$1$};
	\draw (217.2,90.33) node [anchor=north west][inner sep=0.75pt]  [font=\small]  {$\frac{1}{2}$};
	\draw (238.06,155.69) node [anchor=north west][inner sep=0.75pt]  [font=\small]  {$0$};
	\draw (213.99,52.19) node [anchor=north west][inner sep=0.75pt]  [font=\small]  {$1$};
	\draw (340.75,156.49) node [anchor=north west][inner sep=0.75pt]  [font=\small]  {$0$};
	\draw (336.74,26.52) node [anchor=north west][inner sep=0.75pt]  [font=\small]  {$1$};
	\draw (296.64,26.52) node [anchor=north west][inner sep=0.75pt]  [font=\small]  {$CL$};
	\draw (411.02,51.39) node [anchor=north west][inner sep=0.75pt]  [font=\small]  {$\neg ( 1) \ =\ 0$};
	\draw (411.02,73.85) node [anchor=north west][inner sep=0.75pt]  [font=\small]  {$\neg ( 0) \ =\ 1$};
	\draw (403.61,95.25) node [anchor=north west][inner sep=0.75pt]  [font=\small]  {$\neg \left(\frac{1}{2}\right) \ =\ \frac{1}{2}$};
	\draw (187.12,27.13) node [anchor=north west][inner sep=0.75pt]  [font=\normalsize]  {$IL$};
\end{tikzpicture}
\caption{The Intuitionistic and Classical logics}
\end{center}
\end{figure}
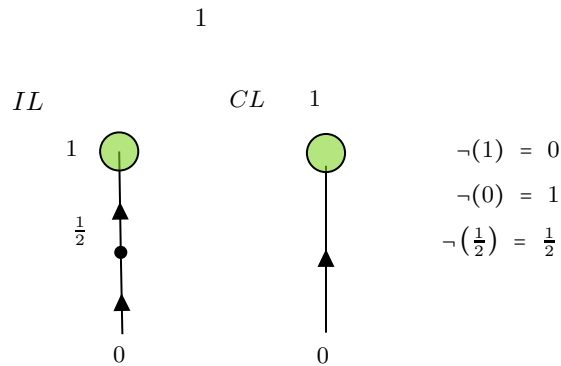

Consider now the following frame:

\begin{figure}[H]
\begin{center}
    \tikzset{every picture/.style={line width=0.75pt}} 
    \begin{tikzpicture}[x=0.75pt,y=0.75pt,yscale=-.7,xscale=.7]
        
        \draw  [fill={rgb, 255:red, 0; green, 0; blue, 0 }  ,fill opacity=1 ] (198,140.5) .. controls (198,138.84) and (199.34,137.5) .. (201,137.5) .. controls (202.66,137.5) and (204,138.84) .. (204,140.5) .. controls (204,142.16) and (202.66,143.5) .. (201,143.5) .. controls (199.34,143.5) and (198,142.16) .. (198,140.5) -- cycle ;
        \draw  [fill={rgb, 255:red, 0; green, 0; blue, 0 }  ,fill opacity=1 ] (338,140.5) .. controls (338,138.84) and (339.34,137.5) .. (341,137.5) .. controls (342.66,137.5) and (344,138.84) .. (344,140.5) .. controls (344,142.16) and (342.66,143.5) .. (341,143.5) .. controls (339.34,143.5) and (338,142.16) .. (338,140.5) -- cycle ;
        \draw    (201,140.5) .. controls (231,148.5) and (217,169.5) .. (202,170.5) .. controls (187.52,171.47) and (184.22,161.26) .. (193,147.06) ;
        \draw [shift={(194,145.5)}, rotate = 483.69] [color={rgb, 255:red, 0; green, 0; blue, 0 }  ][line width=0.75]    (10.93,-3.29) .. controls (6.95,-1.4) and (3.31,-0.3) .. (0,0) .. controls (3.31,0.3) and (6.95,1.4) .. (10.93,3.29)   ;
        \draw    (201,140.5) -- (336,140.5) ;
        \draw [shift={(338,140.5)}, rotate = 180] [color={rgb, 255:red, 0; green, 0; blue, 0 }  ][line width=0.75]    (10.93,-3.29) .. controls (6.95,-1.4) and (3.31,-0.3) .. (0,0) .. controls (3.31,0.3) and (6.95,1.4) .. (10.93,3.29)   ;
        
        \draw (196,109.9) node [anchor=north west][inner sep=0.75pt]    {$w=\langle 1,CL\rangle $};
        \draw (331,110.9) node [anchor=north west][inner sep=0.75pt]    {$w'=\langle 2,IL\rangle $};
    \end{tikzpicture}
\end{center}
\caption{Model with CL and IL}
\end{figure}
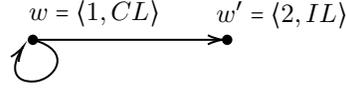

And $v_w(\alpha) = 0$ and $v_{w'}(\alpha) = \frac{1}{2}$, for a propositional variable $\alpha$. In this case we have $w \not \vDash  \Box(\alpha\vee \neg \alpha)$, since  $v_w(\Box(\alpha\vee \neg \alpha))=0$.

\end{ex}

The following example puts together the Logic of Paradox and the Classical Logic.
    
\begin{ex}\

\begin{figure}[H]
\begin{center}
\tikzset{every picture/.style={line width=0.75pt}} 
\begin{tikzpicture}[x=0.75pt,y=0.75pt,yscale=-1,xscale=1]

	\draw    (227.31,181.08) -- (194.3,147.43) ;
	\draw [shift={(210.81,164.26)}, rotate = 45.56] [fill={rgb, 255:red, 0; green, 0; blue, 0 }  ][line width=0.08]  [draw opacity=0] (8.93,-4.29) -- (0,0) -- (8.93,4.29) -- cycle    ;
	\draw    (227.31,181.08) -- (260.32,148.08) ;
	\draw [shift={(243.81,164.58)}, rotate = 135] [fill={rgb, 255:red, 0; green, 0; blue, 0 }  ][line width=0.08]  [draw opacity=0] (8.93,-4.29) -- (0,0) -- (8.93,4.29) -- cycle    ;
	\draw    (260.32,148.08) -- (227.31,105.36) ;
	\draw [shift={(243.81,126.72)}, rotate = 52.31] [fill={rgb, 255:red, 0; green, 0; blue, 0 }  ][line width=0.08]  [draw opacity=0] (8.93,-4.29) -- (0,0) -- (8.93,4.29) -- cycle    ;
	\draw    (194.3,147.43) -- (227.31,105.36) ;
	\draw [shift={(210.81,126.4)}, rotate = 128.12] [fill={rgb, 255:red, 0; green, 0; blue, 0 }  ][line width=0.08]  [draw opacity=0] (8.93,-4.29) -- (0,0) -- (8.93,4.29) -- cycle    ;
	\draw  [fill={rgb, 255:red, 126; green, 211; blue, 33 }  ,fill opacity=0.58 ] (219.54,105.36) .. controls (219.54,101.07) and (223.02,97.6) .. (227.31,97.6) .. controls (231.6,97.6) and (235.08,101.07) .. (235.08,105.36) .. controls (235.08,109.65) and (231.6,113.13) .. (227.31,113.13) .. controls (223.02,113.13) and (219.54,109.65) .. (219.54,105.36) -- cycle ;
	\draw  [fill={rgb, 255:red, 126; green, 211; blue, 33 }  ,fill opacity=0.58 ] (186.54,147.43) .. controls (186.54,143.14) and (190.01,139.66) .. (194.3,139.66) .. controls (198.59,139.66) and (202.07,143.14) .. (202.07,147.43) .. controls (202.07,151.72) and (198.59,155.2) .. (194.3,155.2) .. controls (190.01,155.2) and (186.54,151.72) .. (186.54,147.43) -- cycle ;
	\draw    (335.39,184.32) -- (335.39,116.37) ;
	\draw [shift={(335.39,150.34)}, rotate = 90] [fill={rgb, 255:red, 0; green, 0; blue, 0 }  ][line width=0.08]  [draw opacity=0] (8.93,-4.29) -- (0,0) -- (8.93,4.29) -- cycle    ;
	\draw  [fill={rgb, 255:red, 126; green, 211; blue, 33 }  ,fill opacity=0.58 ] (327.62,108.6) .. controls (327.62,104.31) and (331.1,100.83) .. (335.39,100.83) .. controls (339.68,100.83) and (343.16,104.31) .. (343.16,108.6) .. controls (343.16,112.89) and (339.68,116.37) .. (335.39,116.37) .. controls (331.1,116.37) and (327.62,112.89) .. (327.62,108.6) -- cycle ;

	\draw (145.77,79.47) node [anchor=north west][inner sep=0.75pt]  [font=\normalsize]  {$LP$};
	\draw (143.24,138.71) node [anchor=north west][inner sep=0.75pt]  [font=\small]  {$\{V,F\}$};
	\draw (263.61,138.07) node [anchor=north west][inner sep=0.75pt]  [font=\small]  {$\{\ \}$};
	\draw (216.31,184.02) node [anchor=north west][inner sep=0.75pt]  [font=\small]  {$\{F\}$};
	\draw (215.78,79.82) node [anchor=north west][inner sep=0.75pt]  [font=\small]  {$\{V\}$};
	\draw (324.39,187.25) node [anchor=north west][inner sep=0.75pt]  [font=\small]  {$\{F\}$};
	\draw (323.86,83.05) node [anchor=north west][inner sep=0.75pt]  [font=\small]  {$\{V\}$};
	\draw (279.56,78.17) node [anchor=north west][inner sep=0.75pt]  [font=\normalsize]  {$CL$};
	\draw (379.41,96) node [anchor=north west][inner sep=0.75pt]  [font=\footnotesize]  {$\neg (\{V\}) \ =\ \{F\}$};
	\draw (381.35,114.12) node [anchor=north west][inner sep=0.75pt]  [font=\footnotesize]  {$\neg (\{F\}) \ =\ \{V\}$};
	\draw (386.7,132.24) node [anchor=north west][inner sep=0.75pt]  [font=\footnotesize]  {$\neg (\{\ \}) \ =\ \{\ \}$};
	\draw (361.29,149.07) node [anchor=north west][inner sep=0.75pt]  [font=\footnotesize]  {$\neg (\{V,F\}) \ =\ \{V,F\}$};
\end{tikzpicture}
\caption{The Logic of Paradox and the Classical Logic}
\end{center}
\end{figure}
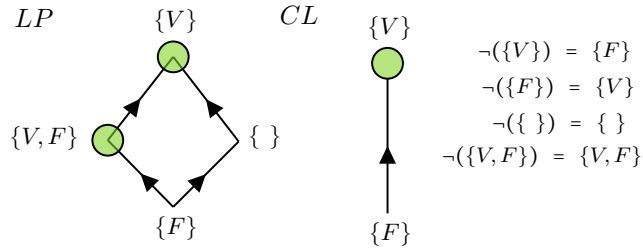

    Take the model

    \begin{figure}[H]
    \begin{center}
    \tikzset{every picture/.style={line width=0.75pt}} 
    \begin{tikzpicture}[x=0.75pt,y=0.75pt,yscale=-1,xscale=1]

        \draw  [fill={rgb, 255:red, 0; green, 0; blue, 0 }  ,fill opacity=1 ] (198,140.5) .. controls (198,138.84) and (199.34,137.5) .. (201,137.5) .. controls (202.66,137.5) and (204,138.84) .. (204,140.5) .. controls (204,142.16) and (202.66,143.5) .. (201,143.5) .. controls (199.34,143.5) and (198,142.16) .. (198,140.5) -- cycle ;
        \draw  [fill={rgb, 255:red, 0; green, 0; blue, 0 }  ,fill opacity=1 ] (338,140.5) .. controls (338,138.84) and (339.34,137.5) .. (341,137.5) .. controls (342.66,137.5) and (344,138.84) .. (344,140.5) .. controls (344,142.16) and (342.66,143.5) .. (341,143.5) .. controls (339.34,143.5) and (338,142.16) .. (338,140.5) -- cycle ;
        \draw    (201,140.5) -- (336,140.5) ;
        \draw [shift={(338,140.5)}, rotate = 180] [color={rgb, 255:red, 0; green, 0; blue, 0 }  ][line width=0.75]    (10.93,-3.29) .. controls (6.95,-1.4) and (3.31,-0.3) .. (0,0) .. controls (3.31,0.3) and (6.95,1.4) .. (10.93,3.29)   ;

        \draw (196,109.9) node [anchor=north west][inner sep=0.75pt]    {$w=\langle 1,CL \rangle $};
        \draw (331,110.9) node [anchor=north west][inner sep=0.75pt]    {$w'=\langle 2,LP \rangle $};
    \end{tikzpicture}
    \end{center}
    \caption{Model with CL and LP}
    \end{figure}

    In this model we have $v_{w'}(\alpha) = \{V,F\}$, for a given propositional variable $\alpha$. It follows that $\alpha$ and $   \neg\alpha$ hold at $w'$, however $v_w(\Box \alpha) =0$.
    
    Consider now the 2-sublattice $CL'$ of $LP$ consisting in the nodes $\{V,F\}$ and $\{F\}$ and the model: 

    \begin{figure}[H]
    \begin{center}
    \tikzset{every picture/.style={line width=0.75pt}} 
    \begin{tikzpicture}[x=0.75pt,y=0.75pt,yscale=-1,xscale=1]

        \draw  [fill={rgb, 255:red, 0; green, 0; blue, 0 }  ,fill opacity=1 ] (198,140.5) .. controls (198,138.84) and (199.34,137.5) .. (201,137.5) .. controls (202.66,137.5) and (204,138.84) .. (204,140.5) .. controls (204,142.16) and (202.66,143.5) .. (201,143.5) .. controls (199.34,143.5) and (198,142.16) .. (198,140.5) -- cycle ;
        \draw  [fill={rgb, 255:red, 0; green, 0; blue, 0 }  ,fill opacity=1 ] (338,140.5) .. controls (338,138.84) and (339.34,137.5) .. (341,137.5) .. controls (342.66,137.5) and (344,138.84) .. (344,140.5) .. controls (344,142.16) and (342.66,143.5) .. (341,143.5) .. controls (339.34,143.5) and (338,142.16) .. (338,140.5) -- cycle ;
        \draw    (201,140.5) -- (336,140.5) ;
        \draw [shift={(338,140.5)}, rotate = 180] [color={rgb, 255:red, 0; green, 0; blue, 0 }  ][line width=0.75]    (10.93,-3.29) .. controls (6.95,-1.4) and (3.31,-0.3) .. (0,0) .. controls (3.31,0.3) and (6.95,1.4) .. (10.93,3.29)   ;

        \draw (196,109.9) node [anchor=north west][inner sep=0.75pt]    {$w=\langle 1,CL' \rangle $};
        \draw (331,110.9) node [anchor=north west][inner sep=0.75pt]    {$w'=\langle 2,LP \rangle $};
    \end{tikzpicture}
    \end{center}
    \caption{Model with CL' and LP}
    \end{figure}

    In this model, $v_{w'}(\alpha) = \{V,F\}$, for a propositional variable $\alpha$. Thus we have that $v_w(\Box (\alpha\wedge \neg \alpha)) =1$.
\end{ex}

Lastly, we will explain how it can be relevant to consider a non-rigid interpretation for negation.

\begin{ex} Consider the lattice $L$ of the extended logic of paradox:

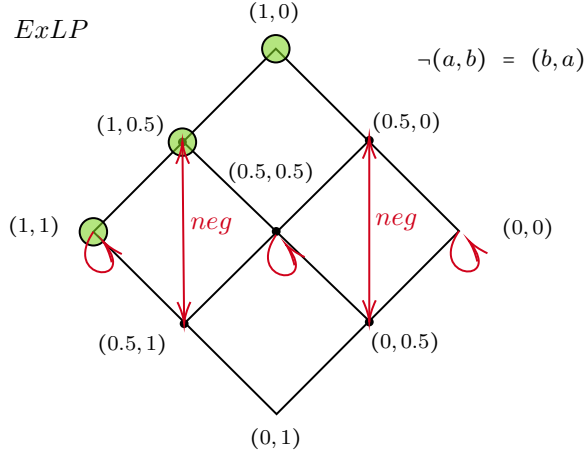
\begin{figure}[H]
\begin{center}
\tikzset{every picture/.style={line width=0.75pt}} 
\begin{tikzpicture}[x=0.75pt,y=0.75pt,yscale=-1,xscale=1]

    \draw   (286.29,30.37) -- (378.68,122.38) -- (286.67,214.78) -- (194.28,122.77) -- cycle ;
    \draw  [fill={rgb, 255:red, 0; green, 0; blue, 0 }  ,fill opacity=1 ] (238.32,169.13) .. controls (238.32,168.17) and (239.09,167.4) .. (240.05,167.4) .. controls (241,167.4) and (241.78,168.17) .. (241.78,169.13) .. controls (241.78,170.08) and (241,170.86) .. (240.05,170.86) .. controls (239.09,170.86) and (238.32,170.08) .. (238.32,169.13) -- cycle ;
    \draw  [fill={rgb, 255:red, 0; green, 0; blue, 0 }  ,fill opacity=1 ] (331.69,168.27) .. controls (331.69,167.31) and (332.47,166.54) .. (333.42,166.54) .. controls (334.38,166.54) and (335.15,167.31) .. (335.15,168.27) .. controls (335.15,169.22) and (334.38,169.99) .. (333.42,169.99) .. controls (332.47,169.99) and (331.69,169.22) .. (331.69,168.27) -- cycle ;
    \draw  [fill={rgb, 255:red, 0; green, 0; blue, 0 }  ,fill opacity=1 ] (331.69,76.62) .. controls (331.69,75.67) and (332.47,74.89) .. (333.42,74.89) .. controls (334.38,74.89) and (335.15,75.67) .. (335.15,76.62) .. controls (335.15,77.58) and (334.38,78.35) .. (333.42,78.35) .. controls (332.47,78.35) and (331.69,77.58) .. (331.69,76.62) -- cycle ;
    \draw  [fill={rgb, 255:red, 0; green, 0; blue, 0 }  ,fill opacity=1 ] (237.46,77.49) .. controls (237.46,76.53) and (238.23,75.76) .. (239.18,75.76) .. controls (240.14,75.76) and (240.91,76.53) .. (240.91,77.49) .. controls (240.91,78.44) and (240.14,79.22) .. (239.18,79.22) .. controls (238.23,79.22) and (237.46,78.44) .. (237.46,77.49) -- cycle ;
    \draw  [fill={rgb, 255:red, 0; green, 0; blue, 0 }  ,fill opacity=1 ] (284.75,122.57) .. controls (284.75,121.62) and (285.53,120.84) .. (286.48,120.84) .. controls (287.44,120.84) and (288.21,121.62) .. (288.21,122.57) .. controls (288.21,123.53) and (287.44,124.3) .. (286.48,124.3) .. controls (285.53,124.3) and (284.75,123.53) .. (284.75,122.57) -- cycle ;
    \draw    (239.36,77.18) -- (333.6,167.96) ;
    \draw    (240.91,168.27) -- (334.29,75.76) ;
    \draw  [fill={rgb, 255:red, 126; green, 211; blue, 33 }  ,fill opacity=0.58 ] (279.37,30.37) .. controls (279.37,26.55) and (282.47,23.45) .. (286.29,23.45) .. controls (290.11,23.45) and (293.21,26.55) .. (293.21,30.37) .. controls (293.21,34.19) and (290.11,37.29) .. (286.29,37.29) .. controls (282.47,37.29) and (279.37,34.19) .. (279.37,30.37) -- cycle ;
    \draw  [fill={rgb, 255:red, 126; green, 211; blue, 33 }  ,fill opacity=0.58 ] (232.27,77.49) .. controls (232.27,73.67) and (235.36,70.57) .. (239.18,70.57) .. controls (243,70.57) and (246.1,73.67) .. (246.1,77.49) .. controls (246.1,81.31) and (243,84.4) .. (239.18,84.4) .. controls (235.36,84.4) and (232.27,81.31) .. (232.27,77.49) -- cycle ;
    \draw  [fill={rgb, 255:red, 126; green, 211; blue, 33 }  ,fill opacity=0.58 ] (187.36,122.77) .. controls (187.36,118.95) and (190.46,115.85) .. (194.28,115.85) .. controls (198.1,115.85) and (201.2,118.95) .. (201.2,122.77) .. controls (201.2,126.59) and (198.1,129.68) .. (194.28,129.68) .. controls (190.46,129.68) and (187.36,126.59) .. (187.36,122.77) -- cycle ;
    \draw [color={rgb, 255:red, 208; green, 2; blue, 27 }  ,draw opacity=1 ]   (239.2,79.49) -- (240.03,167.13) ;
    \draw [shift={(240.05,169.13)}, rotate = 269.46] [color={rgb, 255:red, 208; green, 2; blue, 27 }  ,draw opacity=1 ][line width=0.75]    (10.93,-3.29) .. controls (6.95,-1.4) and (3.31,-0.3) .. (0,0) .. controls (3.31,0.3) and (6.95,1.4) .. (10.93,3.29)   ;
    \draw [shift={(239.18,77.49)}, rotate = 89.46] [color={rgb, 255:red, 208; green, 2; blue, 27 }  ,draw opacity=1 ][line width=0.75]    (10.93,-3.29) .. controls (6.95,-1.4) and (3.31,-0.3) .. (0,0) .. controls (3.31,0.3) and (6.95,1.4) .. (10.93,3.29)   ;
    \draw [color={rgb, 255:red, 208; green, 2; blue, 27 }  ,draw opacity=1 ]   (333.42,76.89) -- (333.42,166.27) ;
    \draw [shift={(333.42,168.27)}, rotate = 270] [color={rgb, 255:red, 208; green, 2; blue, 27 }  ,draw opacity=1 ][line width=0.75]    (10.93,-3.29) .. controls (6.95,-1.4) and (3.31,-0.3) .. (0,0) .. controls (3.31,0.3) and (6.95,1.4) .. (10.93,3.29)   ;
    \draw [shift={(333.42,74.89)}, rotate = 90] [color={rgb, 255:red, 208; green, 2; blue, 27 }  ,draw opacity=1 ][line width=0.75]    (10.93,-3.29) .. controls (6.95,-1.4) and (3.31,-0.3) .. (0,0) .. controls (3.31,0.3) and (6.95,1.4) .. (10.93,3.29)   ;
    \draw [color={rgb, 255:red, 208; green, 2; blue, 27 }  ,draw opacity=1 ]   (378.68,122.38) .. controls (362.29,146.59) and (400.69,150.04) .. (384.62,128.55) ;
    \draw [shift={(383.56,127.2)}, rotate = 50.84] [color={rgb, 255:red, 208; green, 2; blue, 27 }  ,draw opacity=1 ][line width=0.75]    (10.93,-3.29) .. controls (6.95,-1.4) and (3.31,-0.3) .. (0,0) .. controls (3.31,0.3) and (6.95,1.4) .. (10.93,3.29)   ;
    \draw [color={rgb, 255:red, 208; green, 2; blue, 27 }  ,draw opacity=1 ]   (286.48,124.3) .. controls (270.08,148.51) and (308.48,151.96) .. (292.41,130.47) ;
    \draw [shift={(291.36,129.12)}, rotate = 50.84] [color={rgb, 255:red, 208; green, 2; blue, 27 }  ,draw opacity=1 ][line width=0.75]    (10.93,-3.29) .. controls (6.95,-1.4) and (3.31,-0.3) .. (0,0) .. controls (3.31,0.3) and (6.95,1.4) .. (10.93,3.29)   ;
    \draw [color={rgb, 255:red, 208; green, 2; blue, 27 }  ,draw opacity=1 ]   (194.28,122.77) .. controls (177.88,146.98) and (216.28,150.43) .. (200.21,128.93) ;
    \draw [shift={(199.16,127.59)}, rotate = 50.84] [color={rgb, 255:red, 208; green, 2; blue, 27 }  ,draw opacity=1 ][line width=0.75]    (10.93,-3.29) .. controls (6.95,-1.4) and (3.31,-0.3) .. (0,0) .. controls (3.31,0.3) and (6.95,1.4) .. (10.93,3.29)   ;

    \draw (270.8,5.25) node [anchor=north west][inner sep=0.75pt]  [font=\footnotesize]  {$( 1,0)$};
    \draw (260.36,86.09) node [anchor=north west][inner sep=0.75pt]  [font=\footnotesize]  {$( 0.5,0.5)$};
    \draw (333.24,60.58) node [anchor=north west][inner sep=0.75pt]  [font=\footnotesize]  {$( 0.5,0)$};
    \draw (194.04,62.31) node [anchor=north west][inner sep=0.75pt]  [font=\footnotesize]  {$( 1,0.5)$};
    \draw (398.76,114.18) node [anchor=north west][inner sep=0.75pt]  [font=\footnotesize]  {$( 0,0)$};
    \draw (149.77,113.32) node [anchor=north west][inner sep=0.75pt]  [font=\footnotesize]  {$( 1,1)$};
    \draw (332.37,172.11) node [anchor=north west][inner sep=0.75pt]  [font=\footnotesize]  {$( 0,0.5)$};
    \draw (194.91,172.97) node [anchor=north west][inner sep=0.75pt]  [font=\footnotesize]  {$( 0.5,1)$};
    \draw (271.67,220.52) node [anchor=north west][inner sep=0.75pt]  [font=\footnotesize]  {$( 0,1)$};
    \draw (152.55,14.49) node [anchor=north west][inner sep=0.75pt]    {$ExLP$};
    \draw (355.98,28.46) node [anchor=north west][inner sep=0.75pt]  [font=\small]  {$\neg ( a,b) \ =\ ( b,a)$};
    \draw (241.68,113.91) node [anchor=north west][inner sep=0.75pt]  [color={rgb, 255:red, 208; green, 2; blue, 27 }  ,opacity=1 ]  {$neg$};
    \draw (335.05,112.18) node [anchor=north west][inner sep=0.75pt]  [color={rgb, 255:red, 208; green, 2; blue, 27 }  ,opacity=1 ]  {$neg$};
\end{tikzpicture}
\end{center}
\caption{Extended Logic of Paradox}
\end{figure}

If we take a rigid interpretation of the negation, the following lattice in \cref{linear-sublattice} is not a sublattice of $L$. However, if we assume the down interpretation of the negation, the induced negation is the one represented in the lattice, in this way the lattice presented above is, in fact, a sublattice of the extended logic of paradox lattice. A less strict notion of lattice produces more variety of sublattices and, consequently, more interesting many-logic modal structures.  

\begin{figure}[H]
\begin{center}
\tikzset{every picture/.style={line width=0.75pt}} 
\begin{tikzpicture}[x=0.75pt,y=0.75pt,yscale=-1,xscale=1]

    \draw  [fill={rgb, 255:red, 0; green, 0; blue, 0 }  ,fill opacity=1 ] (268.39,126.52) .. controls (268.39,125.91) and (268.88,125.43) .. (269.48,125.43) .. controls (270.08,125.43) and (270.57,125.91) .. (270.57,126.52) .. controls (270.57,127.12) and (270.08,127.61) .. (269.48,127.61) .. controls (268.88,127.61) and (268.39,127.12) .. (268.39,126.52) -- cycle ;
    \draw  [fill={rgb, 255:red, 0; green, 0; blue, 0 }  ,fill opacity=1 ] (327.36,68.09) .. controls (327.36,67.49) and (327.85,67) .. (328.45,67) .. controls (329.05,67) and (329.54,67.49) .. (329.54,68.09) .. controls (329.54,68.69) and (329.05,69.18) .. (328.45,69.18) .. controls (327.85,69.18) and (327.36,68.69) .. (327.36,68.09) -- cycle ;
    \draw    (269.48,127.61) -- (298.69,98.67) -- (328.45,69.18) ;
    \draw  [fill={rgb, 255:red, 126; green, 211; blue, 33 }  ,fill opacity=0.58 ] (294.31,38.88) .. controls (294.31,36.47) and (296.27,34.51) .. (298.68,34.51) .. controls (301.1,34.51) and (303.05,36.47) .. (303.05,38.88) .. controls (303.05,41.29) and (301.1,43.25) .. (298.68,43.25) .. controls (296.27,43.25) and (294.31,41.29) .. (294.31,38.88) -- cycle ;
    \draw [color={rgb, 255:red, 208; green, 2; blue, 27 }  ,draw opacity=1 ]   (328.45,67) -- (299.56,153.45) ;
    \draw [shift={(298.93,155.35)}, rotate = 288.48] [color={rgb, 255:red, 208; green, 2; blue, 27 }  ,draw opacity=1 ][line width=0.75]    (10.93,-3.29) .. controls (6.95,-1.4) and (3.31,-0.3) .. (0,0) .. controls (3.31,0.3) and (6.95,1.4) .. (10.93,3.29)   ;
    \draw    (298.93,155.35) -- (270.88,126.85) ;
    \draw [shift={(269.48,125.43)}, rotate = 45.46] [color={rgb, 255:red, 0; green, 0; blue, 0 }  ][line width=0.75]    (10.93,-3.29) .. controls (6.95,-1.4) and (3.31,-0.3) .. (0,0) .. controls (3.31,0.3) and (6.95,1.4) .. (10.93,3.29)   ;
    \draw    (328.13,68.8) -- (300.09,40.3) ;
    \draw [shift={(298.68,38.88)}, rotate = 45.46] [color={rgb, 255:red, 0; green, 0; blue, 0 }  ][line width=0.75]    (10.93,-3.29) .. controls (6.95,-1.4) and (3.31,-0.3) .. (0,0) .. controls (3.31,0.3) and (6.95,1.4) .. (10.93,3.29)   ;
    \draw [color={rgb, 255:red, 208; green, 2; blue, 27 }  ,draw opacity=1 ]   (268.39,126.52) .. controls (256.37,157.06) and (267.3,149.5) .. (271.13,131.28) ;
    \draw [shift={(271.47,129.56)}, rotate = 100.12] [color={rgb, 255:red, 208; green, 2; blue, 27 }  ,draw opacity=1 ][line width=0.75]    (10.93,-3.29) .. controls (6.95,-1.4) and (3.31,-0.3) .. (0,0) .. controls (3.31,0.3) and (6.95,1.4) .. (10.93,3.29)   ;

    \draw (282.78,20.56) node [anchor=north west][inner sep=0.75pt]  [font=\scriptsize]  {$( 1,0)$};
    \draw (333.27,56.51) node [anchor=north west][inner sep=0.75pt]  [font=\scriptsize]  {$( 0.5,0)$};
    \draw (283.33,156.53) node [anchor=north west][inner sep=0.75pt]  [font=\scriptsize]  {$( 0,1)$};
    \draw (313.59,108.86) node [anchor=north west][inner sep=0.75pt]  [font=\footnotesize,color={rgb, 255:red, 208; green, 2; blue, 27 }  ,opacity=1 ]  {$neg$};
    \draw (238.23,126.88) node [anchor=north west][inner sep=0.75pt]  [font=\footnotesize,color={rgb, 255:red, 208; green, 2; blue, 27 }  ,opacity=1 ]  {$neg$};
\end{tikzpicture}
\end{center}
\caption{Linear sublattice with non-rigid interpretation of negation}
\label{linear-sublattice}
\end{figure}
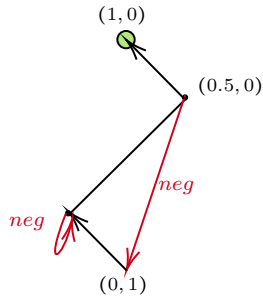

Moreover, once we have down-interpretation of negation, lattices $A_1$ and $A_2$ of \cref{lattice-contain-up-down-lattices} become possible sublattices of the extended Logic of Paradox lattice. The odd appearance of the base lattice in \cref{lattice-contain-up-down-lattices} results from the combined requirement that (i) $A_1$ and $A_2$ are sublattices, and (ii) the base lattice has rigid negation.

\end{ex}

\section{Frames and Graphs}

Traditional modal theory treats the notion of validity in four levels -- namely, (i) world, (ii) model, (iii) frame, (iv) class of frames. After defining satisfaction in a world, one will say that a model satisfies a formula whenever the formula is satisfied in every world of the model. Satisfaction in a frame, however, is obtained differently. It represents truth independently of particular choices in the worlds of a model. This is formalized by removing the valuations from the definition of a model. If a model $M$ is a tuple $\pair{W, R, s}$ where $W$ contains worlds, $R$ is a relation of accessibility between worlds and $s$ is the valuation of atomic formulas in each world in $W$, the frame $F$ of $M$ will be formalized as $\pair{W, R}$. In this case, we say that $F$ satisfy a formula $\varphi$ when all models with a frame $F$ satisfy $\varphi$. 

In this context, classes of frames provide semantics for some important formulas in a modal system. A class of frames $C$ in which all frames have a transitive accessibility relation is such that all frames in $C$ satisfy the formula $\Box \psi \so \Box \Box \psi$. Conversely, if a frame satisfy $\Box \psi \so \Box \Box \psi$, then this is a frame in $C$. Modal formulas, that often represent important metaphysical or epistemological notions, can be understood through structural properties of accessibility between possible worlds. This is a very successful story for a definition. One that a person interested in Kripke semantics is familiar with.
 
But with the addition of a new layer into our modal analysis, we should carefully reevaluate those four levels of validity. Our goal in this section is to study layers of validity once we agreed on the base lattice $L$ where transfer of information between logics is established.
Note that the worlds in our modal structure are now associated to a particular lattice $L_w$, meaning it is a pair composed of an `identifier' for the world and a reference to the lattice in which this world should operate. Our first question is now whether or not we should refer to the lattice $L_w$ in our notion of frame.

When making a decision in this matter, we should take into account two possibilities. First, that a frame should deal with the general notion of modal validity regarding its structural features -- i.e. treating validity with respect only to the arrangement of worlds and the accessibility relation. Second, that our new definition should be a generalization of the traditional one. We choose the second option while adding the notion of graph validity corresponding to the first notion. 

If our notion of frame generalizes the traditional one, we should establish that traditional frames are particular choices in the new framework. This is not possible however if we do not fix the lattices in the worlds of frame as the same. A non-trivial base lattice $L$ will allow for many complete sublattices. So for such $L$ the traditional frames will not be definable. In this case, we should allow that our frames have conditions not only over the accessibility relations but also over what are the lattices in the worlds of our model. Thus we define:

\begin{definition}
For a lattice $L$, we say that a $L$-frame $F$ is the pair $\pair{W, R}$ where $R \subseteq W \times W$ and each $w_{id} \in W$ is a pairs $\pair{id, L_{id}}$ with $L_{id}$ a complete sublattice of $L$.

We say that the $L$-model $M = \pair{W, R, s}$ is a model with frame $F = \{W, R\}$.

For a formula $\varphi$, we say that $F \vDash \varphi$ when for all $L$-modals $M$ with frame $F$ are such that $M \vDash \varphi$.
 \end{definition}

\begin{definition}
If $\mathcal{F}$ is a class of $L$-frames and $\varphi$ a formula, then $\mathcal{F} \vDash \varphi$ when, for all $F \in \mathcal{F}$, we have $F \vDash \varphi$.
\end{definition}

In this context, we can now have a class of $L$-frames that recover the traditional notion of frames, i.e. those in which, for a fixed sublattice $L^*$ of $L$, the pairs $\pair{id, L_{id}}$ are such that $L_{id} = L^*$.

From this definition we can now study interesting classes of frames. We introduce this idea in relation to the current literature on twist structures (introduced independently in \cite{vakarelov1977notes, fidel1977algebraic}; for recent developments, see \cite{Carnielli_Coniglio_2021, fideltwist,ono2014modal, odintsov2009axiomatizing}). A twist lattice $T$ in a broad sense is a lattice obtained by pairs of values ($T \subseteq L \times L$) in a given lattice $L$ with negation. The idea is to incorporate the notion of `affirmative' and `negative' force of a proposition. A value $(a,b) \in T$ for a formula $\varphi$ stands to the idea that we assert $\varphi$ with force $a$ and deny $\varphi$ with force $b$. This concept can be used to produce semantics for families of important logical systems (e.g. Nelson logic introduced in \cite{nelson1949constructible} and systems alike). Let us however work on a narrower framework where the base lattice is a boolean algebra (see. Walter-Coniglio's work on paraconsistent set theory using twist structures built from boolean algebras \cite{Carnielli_Coniglio_2021}).

Consider we have a Boolean algebra $B$. The base twist structure $T$ over $B$ is the set of pairs $\pair{a, b}$ with $a, b \in B$. Conjunction and disjunction are defined by  
$\pair{a,b} \land \pair{c,d} \equiv \pair{a+c, b.d}$ and  
$\pair{a,b} \lor \pair{c,d} \equiv \pair{a.c, b+d}$. The unary function corresponding to negation is the inversion of positive and negative values $\lnot \pair{a, b} \equiv \pair{b, a}$ and the order relation is obtained by $\pair{a, b} \leq \pair{c, d} \equiv (a \leq c \text{ and } b \geq d)$. One may now take lattice subsets of $T$. Notably, the subset $C = \{\pair{a, b} \in T \mid a+b = 1 \text{ and } a.b = 0\}$ corresponds to a lattice equivalent to the original $B$ while the subset $P = \{\pair{a, b} \in T \mid a+b = 1 \}$ corresponds to a `paraconsistent' algebra.

Once we fix a twist lattice structure $T$, we may define a relation of being `more classical' (there are alternative definitions) between its subsets. We say that the pair $\pair{a, b}$ has more excluded middle than $\pair{c, d}$ (notation: $\pair{a, b} \geq^{ExM} \pair{c,d}$) when $a+b \geq c+d$; and we say that $\pair{a, b}$ has more non-contradiction than $\pair{c, d}$ (notation: $\pair{a, b} \geq^{NonC} \pair{c,d}$) when $a.b \leq c+d$. Now, let $T_1$ and $T_2$ be subsets of $T$, we say $T_1$ is `more classical' than $T_2$ (notation: $T_1 \geq ^{Cl} T_2$) when, for all $\pair{x,y} \in T_1$, there is $\pair{z,w} \in T_2$ such that $\pair{x,y} \geq^{ExM} \pair{z,w}$ and $\pair{x,y} \geq^{NonC} \pair{z,w}$. 

Let us consider the following example for the base twist $T$ obtained from the Boolean algebra $L$ in \cref{four-valued-base-lattice}. In the following \cref{base-twist} we add small circles in green next to $T$'s values for the sublattice $T_B$ (equivalent to the original Boolean algebra), in red for $T_{(a)}$ and in red for $T_{(1)}$:

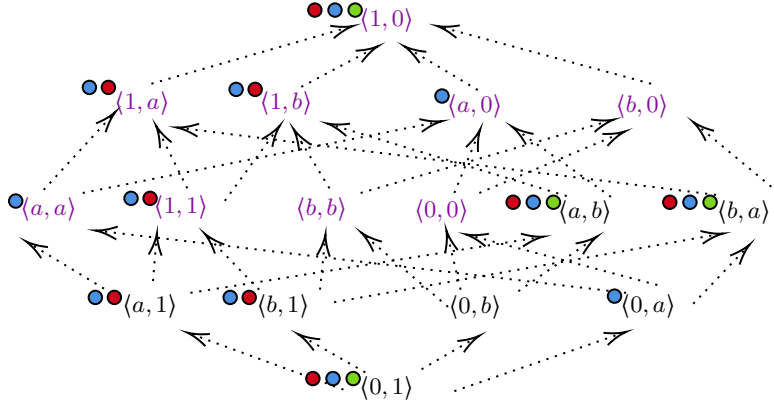
\begin{figure}[H]
\begin{center}
\tikzset{every picture/.style={line width=0.75pt}} 
\begin{tikzpicture}[x=0.75pt,y=0.75pt,yscale=-1,xscale=1]

    \draw  [dash pattern={on 0.84pt off 2.51pt}]  (218.07,265.28) -- (141.2,237.22) ;
    \draw [shift={(139.32,236.53)}, rotate = 20.05] [color={rgb, 255:red, 0; green, 0; blue, 0 }  ][line width=0.75]    (10.93,-3.29) .. controls (6.95,-1.4) and (3.31,-0.3) .. (0,0) .. controls (3.31,0.3) and (6.95,1.4) .. (10.93,3.29)   ;
    \draw  [dash pattern={on 0.84pt off 2.51pt}]  (273.99,266.07) -- (358.73,236.79) ;
    \draw [shift={(360.62,236.14)}, rotate = 160.94] [color={rgb, 255:red, 0; green, 0; blue, 0 }  ][line width=0.75]    (10.93,-3.29) .. controls (6.95,-1.4) and (3.31,-0.3) .. (0,0) .. controls (3.31,0.3) and (6.95,1.4) .. (10.93,3.29)   ;
    \draw  [dash pattern={on 0.84pt off 2.51pt}]  (229.1,256.62) -- (196.97,238.31) ;
    \draw [shift={(195.23,237.32)}, rotate = 29.67] [color={rgb, 255:red, 0; green, 0; blue, 0 }  ][line width=0.75]    (10.93,-3.29) .. controls (6.95,-1.4) and (3.31,-0.3) .. (0,0) .. controls (3.31,0.3) and (6.95,1.4) .. (10.93,3.29)   ;
    \draw  [dash pattern={on 0.84pt off 2.51pt}]  (257.45,255.04) -- (284.86,239.11) ;
    \draw [shift={(286.59,238.11)}, rotate = 149.84] [color={rgb, 255:red, 0; green, 0; blue, 0 }  ][line width=0.75]    (10.93,-3.29) .. controls (6.95,-1.4) and (3.31,-0.3) .. (0,0) .. controls (3.31,0.3) and (6.95,1.4) .. (10.93,3.29)   ;
    \draw  [dash pattern={on 0.84pt off 2.51pt}]  (122,210.5) -- (123.83,188.49) ;
    \draw [shift={(124,186.5)}, rotate = 94.76] [color={rgb, 255:red, 0; green, 0; blue, 0 }  ][line width=0.75]    (10.93,-3.29) .. controls (6.95,-1.4) and (3.31,-0.3) .. (0,0) .. controls (3.31,0.3) and (6.95,1.4) .. (10.93,3.29)   ;
    \draw  [dash pattern={on 0.84pt off 2.51pt}]  (101.14,216.66) -- (59.71,191.54) ;
    \draw [shift={(58,190.5)}, rotate = 31.24] [color={rgb, 255:red, 0; green, 0; blue, 0 }  ][line width=0.75]    (10.93,-3.29) .. controls (6.95,-1.4) and (3.31,-0.3) .. (0,0) .. controls (3.31,0.3) and (6.95,1.4) .. (10.93,3.29)   ;
    \draw  [dash pattern={on 0.84pt off 2.51pt}]  (140.89,216.45) -- (316.9,188.42) ;
    \draw [shift={(318.88,188.1)}, rotate = 170.95] [color={rgb, 255:red, 0; green, 0; blue, 0 }  ][line width=0.75]    (10.93,-3.29) .. controls (6.95,-1.4) and (3.31,-0.3) .. (0,0) .. controls (3.31,0.3) and (6.95,1.4) .. (10.93,3.29)   ;
    \draw  [dash pattern={on 0.84pt off 2.51pt}]  (175.55,214.88) -- (149.39,187.94) ;
    \draw [shift={(148,186.5)}, rotate = 45.85] [color={rgb, 255:red, 0; green, 0; blue, 0 }  ][line width=0.75]    (10.93,-3.29) .. controls (6.95,-1.4) and (3.31,-0.3) .. (0,0) .. controls (3.31,0.3) and (6.95,1.4) .. (10.93,3.29)   ;
    \draw  [dash pattern={on 0.84pt off 2.51pt}]  (206.26,214.88) -- (208.8,188.99) ;
    \draw [shift={(209,187)}, rotate = 95.61] [color={rgb, 255:red, 0; green, 0; blue, 0 }  ][line width=0.75]    (10.93,-3.29) .. controls (6.95,-1.4) and (3.31,-0.3) .. (0,0) .. controls (3.31,0.3) and (6.95,1.4) .. (10.93,3.29)   ;
    \draw  [dash pattern={on 0.84pt off 2.51pt}]  (213.35,221.18) -- (407.48,186.87) ;
    \draw [shift={(409.45,186.53)}, rotate = 169.98] [color={rgb, 255:red, 0; green, 0; blue, 0 }  ][line width=0.75]    (10.93,-3.29) .. controls (6.95,-1.4) and (3.31,-0.3) .. (0,0) .. controls (3.31,0.3) and (6.95,1.4) .. (10.93,3.29)   ;
    \draw  [dash pattern={on 0.84pt off 2.51pt}]  (270.84,223.54) -- (228.52,187.3) ;
    \draw [shift={(227,186)}, rotate = 40.57] [color={rgb, 255:red, 0; green, 0; blue, 0 }  ][line width=0.75]    (10.93,-3.29) .. controls (6.95,-1.4) and (3.31,-0.3) .. (0,0) .. controls (3.31,0.3) and (6.95,1.4) .. (10.93,3.29)   ;
    \draw  [dash pattern={on 0.84pt off 2.51pt}]  (277.14,214.88) -- (270.46,186.95) ;
    \draw [shift={(270,185)}, rotate = 76.56] [color={rgb, 255:red, 0; green, 0; blue, 0 }  ][line width=0.75]    (10.93,-3.29) .. controls (6.95,-1.4) and (3.31,-0.3) .. (0,0) .. controls (3.31,0.3) and (6.95,1.4) .. (10.93,3.29)   ;
    \draw  [dash pattern={on 0.84pt off 2.51pt}]  (306.28,215.66) -- (344.36,189.14) ;
    \draw [shift={(346,188)}, rotate = 145.14] [color={rgb, 255:red, 0; green, 0; blue, 0 }  ][line width=0.75]    (10.93,-3.29) .. controls (6.95,-1.4) and (3.31,-0.3) .. (0,0) .. controls (3.31,0.3) and (6.95,1.4) .. (10.93,3.29)   ;
    \draw  [dash pattern={on 0.84pt off 2.51pt}]  (360.62,216.45) -- (94.99,185.23) ;
    \draw [shift={(93,185)}, rotate = 6.7] [color={rgb, 255:red, 0; green, 0; blue, 0 }  ][line width=0.75]    (10.93,-3.29) .. controls (6.95,-1.4) and (3.31,-0.3) .. (0,0) .. controls (3.31,0.3) and (6.95,1.4) .. (10.93,3.29)   ;
    \draw  [dash pattern={on 0.84pt off 2.51pt}]  (376.37,212.51) -- (281.93,187.51) ;
    \draw [shift={(280,187)}, rotate = 14.83] [color={rgb, 255:red, 0; green, 0; blue, 0 }  ][line width=0.75]    (10.93,-3.29) .. controls (6.95,-1.4) and (3.31,-0.3) .. (0,0) .. controls (3.31,0.3) and (6.95,1.4) .. (10.93,3.29)   ;
    \draw  [dash pattern={on 0.84pt off 2.51pt}]  (394.48,221.97) -- (421.62,193.45) ;
    \draw [shift={(423,192)}, rotate = 133.58] [color={rgb, 255:red, 0; green, 0; blue, 0 }  ][line width=0.75]    (10.93,-3.29) .. controls (6.95,-1.4) and (3.31,-0.3) .. (0,0) .. controls (3.31,0.3) and (6.95,1.4) .. (10.93,3.29)   ;
    \draw  [dash pattern={on 0.84pt off 2.51pt}]  (66.86,164.47) -- (101.74,126.57) ;
    \draw [shift={(103.09,125.1)}, rotate = 132.61] [color={rgb, 255:red, 0; green, 0; blue, 0 }  ][line width=0.75]    (10.93,-3.29) .. controls (6.95,-1.4) and (3.31,-0.3) .. (0,0) .. controls (3.31,0.3) and (6.95,1.4) .. (10.93,3.29)   ;
    \draw  [dash pattern={on 0.84pt off 2.51pt}]  (86.55,166.05) -- (265.73,130.21) ;
    \draw [shift={(267.69,129.82)}, rotate = 168.69] [color={rgb, 255:red, 0; green, 0; blue, 0 }  ][line width=0.75]    (10.93,-3.29) .. controls (6.95,-1.4) and (3.31,-0.3) .. (0,0) .. controls (3.31,0.3) and (6.95,1.4) .. (10.93,3.29)   ;
    \draw  [dash pattern={on 0.84pt off 2.51pt}]  (139.32,163.69) -- (124.46,133.97) ;
    \draw [shift={(123.57,132.18)}, rotate = 63.43] [color={rgb, 255:red, 0; green, 0; blue, 0 }  ][line width=0.75]    (10.93,-3.29) .. controls (6.95,-1.4) and (3.31,-0.3) .. (0,0) .. controls (3.31,0.3) and (6.95,1.4) .. (10.93,3.29)   ;
    \draw  [dash pattern={on 0.84pt off 2.51pt}]  (158.22,169.99) -- (184.82,133.61) ;
    \draw [shift={(186,132)}, rotate = 126.18] [color={rgb, 255:red, 0; green, 0; blue, 0 }  ][line width=0.75]    (10.93,-3.29) .. controls (6.95,-1.4) and (3.31,-0.3) .. (0,0) .. controls (3.31,0.3) and (6.95,1.4) .. (10.93,3.29)   ;
    \draw  [dash pattern={on 0.84pt off 2.51pt}]  (210.99,163.69) -- (193.91,135.47) ;
    \draw [shift={(192.87,133.76)}, rotate = 58.82] [color={rgb, 255:red, 0; green, 0; blue, 0 }  ][line width=0.75]    (10.93,-3.29) .. controls (6.95,-1.4) and (3.31,-0.3) .. (0,0) .. controls (3.31,0.3) and (6.95,1.4) .. (10.93,3.29)   ;
    \draw  [dash pattern={on 0.84pt off 2.51pt}]  (226.74,166.05) -- (352.4,129.59) ;
    \draw [shift={(354.32,129.03)}, rotate = 163.82] [color={rgb, 255:red, 0; green, 0; blue, 0 }  ][line width=0.75]    (10.93,-3.29) .. controls (6.95,-1.4) and (3.31,-0.3) .. (0,0) .. controls (3.31,0.3) and (6.95,1.4) .. (10.93,3.29)   ;
    \draw  [dash pattern={on 0.84pt off 2.51pt}]  (273.99,166.05) -- (286.25,135.85) ;
    \draw [shift={(287,134)}, rotate = 112.1] [color={rgb, 255:red, 0; green, 0; blue, 0 }  ][line width=0.75]    (10.93,-3.29) .. controls (6.95,-1.4) and (3.31,-0.3) .. (0,0) .. controls (3.31,0.3) and (6.95,1.4) .. (10.93,3.29)   ;
    \draw  [dash pattern={on 0.84pt off 2.51pt}]  (286.59,166.05) -- (358.16,135.78) ;
    \draw [shift={(360,135)}, rotate = 157.07] [color={rgb, 255:red, 0; green, 0; blue, 0 }  ][line width=0.75]    (10.93,-3.29) .. controls (6.95,-1.4) and (3.31,-0.3) .. (0,0) .. controls (3.31,0.3) and (6.95,1.4) .. (10.93,3.29)   ;
    \draw  [dash pattern={on 0.84pt off 2.51pt}]  (330.69,167.62) -- (211.91,131.58) ;
    \draw [shift={(210,131)}, rotate = 16.88] [color={rgb, 255:red, 0; green, 0; blue, 0 }  ][line width=0.75]    (10.93,-3.29) .. controls (6.95,-1.4) and (3.31,-0.3) .. (0,0) .. controls (3.31,0.3) and (6.95,1.4) .. (10.93,3.29)   ;
    \draw  [dash pattern={on 0.84pt off 2.51pt}]  (353,168) -- (304.76,133.74) ;
    \draw [shift={(303.13,132.58)}, rotate = 35.38] [color={rgb, 255:red, 0; green, 0; blue, 0 }  ][line width=0.75]    (10.93,-3.29) .. controls (6.95,-1.4) and (3.31,-0.3) .. (0,0) .. controls (3.31,0.3) and (6.95,1.4) .. (10.93,3.29)   ;
    \draw  [dash pattern={on 0.84pt off 2.51pt}]  (406,167) -- (135.98,133.25) ;
    \draw [shift={(134,133)}, rotate = 7.13] [color={rgb, 255:red, 0; green, 0; blue, 0 }  ][line width=0.75]    (10.93,-3.29) .. controls (6.95,-1.4) and (3.31,-0.3) .. (0,0) .. controls (3.31,0.3) and (6.95,1.4) .. (10.93,3.29)   ;
    \draw  [dash pattern={on 0.84pt off 2.51pt}]  (433.86,162.9) -- (396.02,131.49) ;
    \draw [shift={(394.48,130.22)}, rotate = 39.69] [color={rgb, 255:red, 0; green, 0; blue, 0 }  ][line width=0.75]    (10.93,-3.29) .. controls (6.95,-1.4) and (3.31,-0.3) .. (0,0) .. controls (3.31,0.3) and (6.95,1.4) .. (10.93,3.29)   ;
    \draw  [dash pattern={on 0.84pt off 2.51pt}]  (121.21,110.92) -- (220.88,82.33) ;
    \draw [shift={(222.8,81.78)}, rotate = 164] [color={rgb, 255:red, 0; green, 0; blue, 0 }  ][line width=0.75]    (10.93,-3.29) .. controls (6.95,-1.4) and (3.31,-0.3) .. (0,0) .. controls (3.31,0.3) and (6.95,1.4) .. (10.93,3.29)   ;
    \draw  [dash pattern={on 0.84pt off 2.51pt}]  (196.81,112.5) -- (233.63,92.96) ;
    \draw [shift={(235.4,92.02)}, rotate = 152.05] [color={rgb, 255:red, 0; green, 0; blue, 0 }  ][line width=0.75]    (10.93,-3.29) .. controls (6.95,-1.4) and (3.31,-0.3) .. (0,0) .. controls (3.31,0.3) and (6.95,1.4) .. (10.93,3.29)   ;
    \draw  [dash pattern={on 0.84pt off 2.51pt}]  (285.8,114.07) -- (250.51,93.04) ;
    \draw [shift={(248.79,92.02)}, rotate = 30.78] [color={rgb, 255:red, 0; green, 0; blue, 0 }  ][line width=0.75]    (10.93,-3.29) .. controls (6.95,-1.4) and (3.31,-0.3) .. (0,0) .. controls (3.31,0.3) and (6.95,1.4) .. (10.93,3.29)   ;
    \draw  [dash pattern={on 0.84pt off 2.51pt}]  (373.22,110.92) -- (267.26,82.3) ;
    \draw [shift={(265.33,81.78)}, rotate = 15.11] [color={rgb, 255:red, 0; green, 0; blue, 0 }  ][line width=0.75]    (10.93,-3.29) .. controls (6.95,-1.4) and (3.31,-0.3) .. (0,0) .. controls (3.31,0.3) and (6.95,1.4) .. (10.93,3.29)   ;
    \draw  [fill={rgb, 255:red, 208; green, 2; blue, 27 }  ,fill opacity=1 ] (96.15,112.76) .. controls (96.15,110.89) and (97.66,109.38) .. (99.53,109.38) .. controls (101.4,109.38) and (102.92,110.89) .. (102.92,112.76) .. controls (102.92,114.63) and (101.4,116.15) .. (99.53,116.15) .. controls (97.66,116.15) and (96.15,114.63) .. (96.15,112.76) -- cycle ;
    \draw  [fill={rgb, 255:red, 74; green, 144; blue, 226 }  ,fill opacity=1 ] (86.55,112.76) .. controls (86.55,110.89) and (88.07,109.38) .. (89.94,109.38) .. controls (91.81,109.38) and (93.32,110.89) .. (93.32,112.76) .. controls (93.32,114.63) and (91.81,116.15) .. (89.94,116.15) .. controls (88.07,116.15) and (86.55,114.63) .. (86.55,112.76) -- cycle ;
    \draw  [fill={rgb, 255:red, 208; green, 2; blue, 27 }  ,fill opacity=1 ] (300.26,170.73) .. controls (300.26,168.86) and (301.77,167.34) .. (303.64,167.34) .. controls (305.51,167.34) and (307.03,168.86) .. (307.03,170.73) .. controls (307.03,172.6) and (305.51,174.11) .. (303.64,174.11) .. controls (301.77,174.11) and (300.26,172.6) .. (300.26,170.73) -- cycle ;
    \draw  [fill={rgb, 255:red, 74; green, 144; blue, 226 }  ,fill opacity=1 ] (310.42,170.73) .. controls (310.42,168.86) and (311.93,167.34) .. (313.8,167.34) .. controls (315.67,167.34) and (317.19,168.86) .. (317.19,170.73) .. controls (317.19,172.6) and (315.67,174.11) .. (313.8,174.11) .. controls (311.93,174.11) and (310.42,172.6) .. (310.42,170.73) -- cycle ;
    \draw  [fill={rgb, 255:red, 126; green, 211; blue, 33 }  ,fill opacity=1 ] (320.58,170.73) .. controls (320.58,168.86) and (322.09,167.34) .. (323.96,167.34) .. controls (325.83,167.34) and (327.35,168.86) .. (327.35,170.73) .. controls (327.35,172.6) and (325.83,174.11) .. (323.96,174.11) .. controls (322.09,174.11) and (320.58,172.6) .. (320.58,170.73) -- cycle ;
    \draw  [fill={rgb, 255:red, 74; green, 144; blue, 226 }  ,fill opacity=1 ] (264.38,117.11) .. controls (264.38,115.24) and (265.9,113.73) .. (267.77,113.73) .. controls (269.64,113.73) and (271.16,115.24) .. (271.16,117.11) .. controls (271.16,118.98) and (269.64,120.5) .. (267.77,120.5) .. controls (265.9,120.5) and (264.38,118.98) .. (264.38,117.11) -- cycle ;
    \draw  [fill={rgb, 255:red, 208; green, 2; blue, 27 }  ,fill opacity=1 ] (117.15,168.76) .. controls (117.15,166.89) and (118.66,165.38) .. (120.53,165.38) .. controls (122.4,165.38) and (123.92,166.89) .. (123.92,168.76) .. controls (123.92,170.63) and (122.4,172.15) .. (120.53,172.15) .. controls (118.66,172.15) and (117.15,170.63) .. (117.15,168.76) -- cycle ;
    \draw  [fill={rgb, 255:red, 74; green, 144; blue, 226 }  ,fill opacity=1 ] (107.55,168.76) .. controls (107.55,166.89) and (109.07,165.38) .. (110.94,165.38) .. controls (112.81,165.38) and (114.32,166.89) .. (114.32,168.76) .. controls (114.32,170.63) and (112.81,172.15) .. (110.94,172.15) .. controls (109.07,172.15) and (107.55,170.63) .. (107.55,168.76) -- cycle ;
    \draw  [fill={rgb, 255:red, 208; green, 2; blue, 27 }  ,fill opacity=1 ] (170.15,113.76) .. controls (170.15,111.89) and (171.66,110.38) .. (173.53,110.38) .. controls (175.4,110.38) and (176.92,111.89) .. (176.92,113.76) .. controls (176.92,115.63) and (175.4,117.15) .. (173.53,117.15) .. controls (171.66,117.15) and (170.15,115.63) .. (170.15,113.76) -- cycle ;
    \draw  [fill={rgb, 255:red, 74; green, 144; blue, 226 }  ,fill opacity=1 ] (160.55,113.76) .. controls (160.55,111.89) and (162.07,110.38) .. (163.94,110.38) .. controls (165.81,110.38) and (167.32,111.89) .. (167.32,113.76) .. controls (167.32,115.63) and (165.81,117.15) .. (163.94,117.15) .. controls (162.07,117.15) and (160.55,115.63) .. (160.55,113.76) -- cycle ;
    \draw  [fill={rgb, 255:red, 208; green, 2; blue, 27 }  ,fill opacity=1 ] (99.15,218.76) .. controls (99.15,216.89) and (100.66,215.38) .. (102.53,215.38) .. controls (104.4,215.38) and (105.92,216.89) .. (105.92,218.76) .. controls (105.92,220.63) and (104.4,222.15) .. (102.53,222.15) .. controls (100.66,222.15) and (99.15,220.63) .. (99.15,218.76) -- cycle ;
    \draw  [fill={rgb, 255:red, 74; green, 144; blue, 226 }  ,fill opacity=1 ] (89.55,218.76) .. controls (89.55,216.89) and (91.07,215.38) .. (92.94,215.38) .. controls (94.81,215.38) and (96.32,216.89) .. (96.32,218.76) .. controls (96.32,220.63) and (94.81,222.15) .. (92.94,222.15) .. controls (91.07,222.15) and (89.55,220.63) .. (89.55,218.76) -- cycle ;
    \draw  [fill={rgb, 255:red, 208; green, 2; blue, 27 }  ,fill opacity=1 ] (167.15,218.76) .. controls (167.15,216.89) and (168.66,215.38) .. (170.53,215.38) .. controls (172.4,215.38) and (173.92,216.89) .. (173.92,218.76) .. controls (173.92,220.63) and (172.4,222.15) .. (170.53,222.15) .. controls (168.66,222.15) and (167.15,220.63) .. (167.15,218.76) -- cycle ;
    \draw  [fill={rgb, 255:red, 74; green, 144; blue, 226 }  ,fill opacity=1 ] (157.55,218.76) .. controls (157.55,216.89) and (159.07,215.38) .. (160.94,215.38) .. controls (162.81,215.38) and (164.32,216.89) .. (164.32,218.76) .. controls (164.32,220.63) and (162.81,222.15) .. (160.94,222.15) .. controls (159.07,222.15) and (157.55,220.63) .. (157.55,218.76) -- cycle ;
    \draw  [fill={rgb, 255:red, 74; green, 144; blue, 226 }  ,fill opacity=1 ] (49.38,170.11) .. controls (49.38,168.24) and (50.9,166.73) .. (52.77,166.73) .. controls (54.64,166.73) and (56.16,168.24) .. (56.16,170.11) .. controls (56.16,171.98) and (54.64,173.5) .. (52.77,173.5) .. controls (50.9,173.5) and (49.38,171.98) .. (49.38,170.11) -- cycle ;
    \draw  [fill={rgb, 255:red, 74; green, 144; blue, 226 }  ,fill opacity=1 ] (351.38,218.11) .. controls (351.38,216.24) and (352.9,214.73) .. (354.77,214.73) .. controls (356.64,214.73) and (358.16,216.24) .. (358.16,218.11) .. controls (358.16,219.98) and (356.64,221.5) .. (354.77,221.5) .. controls (352.9,221.5) and (351.38,219.98) .. (351.38,218.11) -- cycle ;
    \draw  [fill={rgb, 255:red, 208; green, 2; blue, 27 }  ,fill opacity=1 ] (199.26,259.73) .. controls (199.26,257.86) and (200.77,256.34) .. (202.64,256.34) .. controls (204.51,256.34) and (206.03,257.86) .. (206.03,259.73) .. controls (206.03,261.6) and (204.51,263.11) .. (202.64,263.11) .. controls (200.77,263.11) and (199.26,261.6) .. (199.26,259.73) -- cycle ;
    \draw  [fill={rgb, 255:red, 74; green, 144; blue, 226 }  ,fill opacity=1 ] (209.42,259.73) .. controls (209.42,257.86) and (210.93,256.34) .. (212.8,256.34) .. controls (214.67,256.34) and (216.19,257.86) .. (216.19,259.73) .. controls (216.19,261.6) and (214.67,263.11) .. (212.8,263.11) .. controls (210.93,263.11) and (209.42,261.6) .. (209.42,259.73) -- cycle ;
    \draw  [fill={rgb, 255:red, 126; green, 211; blue, 33 }  ,fill opacity=1 ] (219.58,259.73) .. controls (219.58,257.86) and (221.09,256.34) .. (222.96,256.34) .. controls (224.83,256.34) and (226.35,257.86) .. (226.35,259.73) .. controls (226.35,261.6) and (224.83,263.11) .. (222.96,263.11) .. controls (221.09,263.11) and (219.58,261.6) .. (219.58,259.73) -- cycle ;
    \draw  [fill={rgb, 255:red, 208; green, 2; blue, 27 }  ,fill opacity=1 ] (379.26,170.73) .. controls (379.26,168.86) and (380.77,167.34) .. (382.64,167.34) .. controls (384.51,167.34) and (386.03,168.86) .. (386.03,170.73) .. controls (386.03,172.6) and (384.51,174.11) .. (382.64,174.11) .. controls (380.77,174.11) and (379.26,172.6) .. (379.26,170.73) -- cycle ;
    \draw  [fill={rgb, 255:red, 74; green, 144; blue, 226 }  ,fill opacity=1 ] (389.42,170.73) .. controls (389.42,168.86) and (390.93,167.34) .. (392.8,167.34) .. controls (394.67,167.34) and (396.19,168.86) .. (396.19,170.73) .. controls (396.19,172.6) and (394.67,174.11) .. (392.8,174.11) .. controls (390.93,174.11) and (389.42,172.6) .. (389.42,170.73) -- cycle ;
    \draw  [fill={rgb, 255:red, 126; green, 211; blue, 33 }  ,fill opacity=1 ] (399.58,170.73) .. controls (399.58,168.86) and (401.09,167.34) .. (402.96,167.34) .. controls (404.83,167.34) and (406.35,168.86) .. (406.35,170.73) .. controls (406.35,172.6) and (404.83,174.11) .. (402.96,174.11) .. controls (401.09,174.11) and (399.58,172.6) .. (399.58,170.73) -- cycle ;
    \draw  [fill={rgb, 255:red, 208; green, 2; blue, 27 }  ,fill opacity=1 ] (200.26,73.73) .. controls (200.26,71.86) and (201.77,70.34) .. (203.64,70.34) .. controls (205.51,70.34) and (207.03,71.86) .. (207.03,73.73) .. controls (207.03,75.6) and (205.51,77.11) .. (203.64,77.11) .. controls (201.77,77.11) and (200.26,75.6) .. (200.26,73.73) -- cycle ;
    \draw  [fill={rgb, 255:red, 74; green, 144; blue, 226 }  ,fill opacity=1 ] (210.42,73.73) .. controls (210.42,71.86) and (211.93,70.34) .. (213.8,70.34) .. controls (215.67,70.34) and (217.19,71.86) .. (217.19,73.73) .. controls (217.19,75.6) and (215.67,77.11) .. (213.8,77.11) .. controls (211.93,77.11) and (210.42,75.6) .. (210.42,73.73) -- cycle ;
    \draw  [fill={rgb, 255:red, 126; green, 211; blue, 33 }  ,fill opacity=1 ] (220.58,73.73) .. controls (220.58,71.86) and (222.09,70.34) .. (223.96,70.34) .. controls (225.83,70.34) and (227.35,71.86) .. (227.35,73.73) .. controls (227.35,75.6) and (225.83,77.11) .. (223.96,77.11) .. controls (222.09,77.11) and (220.58,75.6) .. (220.58,73.73) -- cycle ;

    \draw (174.55,113.06) node [anchor=north west][inner sep=0.75pt]  [font=\small,color={rgb, 255:red, 134; green, 14; blue, 158 }  ,opacity=1 ]  {$\langle 1,b\rangle $};
    \draw (252.52,167.4) node [anchor=north west][inner sep=0.75pt]  [font=\small,color={rgb, 255:red, 134; green, 14; blue, 158 }  ,opacity=1 ]  {$\langle 0,0\rangle $};
    \draw (121.11,165.82) node [anchor=north west][inner sep=0.75pt]  [font=\small,color={rgb, 255:red, 134; green, 14; blue, 158 }  ,opacity=1 ]  {$\langle 1,1\rangle $};
    \draw (224.85,71.32) node [anchor=north west][inner sep=0.75pt]  [font=\small,color={rgb, 255:red, 134; green, 14; blue, 158 }  ,opacity=1 ]  {$\langle 1,0\rangle $};
    \draw (354.11,114.63) node [anchor=north west][inner sep=0.75pt]  [font=\small,color={rgb, 255:red, 134; green, 14; blue, 158 }  ,opacity=1 ]  {$\langle b,0\rangle $};
    \draw (269.06,114.63) node [anchor=north west][inner sep=0.75pt]  [font=\small,color={rgb, 255:red, 134; green, 14; blue, 158 }  ,opacity=1 ]  {$\langle a,0\rangle $};
    \draw (101.31,113.06) node [anchor=north west][inner sep=0.75pt]  [font=\small,color={rgb, 255:red, 134; green, 14; blue, 158 }  ,opacity=1 ]  {$\langle 1,a\rangle $};
    \draw (192.77,166.61) node [anchor=north west][inner sep=0.75pt]  [font=\small,color={rgb, 255:red, 134; green, 14; blue, 158 }  ,opacity=1 ]  {$\langle b,b\rangle $};
    \draw (325.08,168.19) node [anchor=north west][inner sep=0.75pt]  [font=\small]  {$\langle a,b\rangle $};
    \draw (405.2,168.19) node [anchor=north west][inner sep=0.75pt]  [font=\small]  {$\langle b,a\rangle $};
    \draw (53.9,166.69) node [anchor=north west][inner sep=0.75pt]  [font=\small,color={rgb, 255:red, 134; green, 14; blue, 158 }  ,opacity=1 ]  {$\langle a,a\rangle $};
    \draw (356.48,215.44) node [anchor=north west][inner sep=0.75pt]  [font=\small]  {$\langle 0,a\rangle $};
    \draw (269.85,216.23) node [anchor=north west][inner sep=0.75pt]  [font=\small]  {$\langle 0,b\rangle $};
    \draw (105.25,215.44) node [anchor=north west][inner sep=0.75pt]  [font=\small]  {$\langle a,1\rangle $};
    \draw (172.98,215.44) node [anchor=north west][inner sep=0.75pt]  [font=\small]  {$\langle b,1\rangle $};
    \draw (224.85,257.18) node [anchor=north west][inner sep=0.75pt]  [font=\small]  {$\langle 0,1\rangle $};
\end{tikzpicture}
\caption{Base twist $T$ for the Boolean algebra $L$}
\label{base-twist}
\end{center}
\end{figure}

We use the notation $T_{(z)}$ for the twist subset $\{\pair{x,y} \mid x+y \geq z\}$. Thus $T_{(a)}$ is the set $\{\pair{x,y} \mid x+y \geq a\}$ while $T_{(1)}$ is the set $\{\pair{x,y} \mid x+y \geq 1\}$. It is then easy to verify that $T_B \geq^{Cl} T_{(1)} \geq^{Cl} T_{(a)} \geq^{Cl} T_{(0)}$ and $T_{(1)} \geq^{Cl} T_{(b)} \geq^{Cl} T_{(0)}$ while $T_{(a)}$ and $T_{(b)}$ cannot be compared with $\geq^{Cl}$. With the filter $Tr = \{\pair{1,0}, \pair{1,a}, \pair{1,b}, \pair{a,0}, \pair{b,0}, \pair{a,a}, \pair{1,1}, \pair{b,b}, \pair{0,0}\}$ (purple values in the figure), we define the following classes of frames:

\begin{enumerate}
    \item[(1)] Classically increasing $\mathcal{F}_{inc}$: $F \in \mathcal{F}_{inc}$ if, and only if, worlds in $F$ have lattice $T_{(z)}$ for some $z \in \{1, a, b, 0\}$, the accessibility relation $R$ is transitive\footnote{This is not required for the analysis proposed. We add transitivity to complete the picture of a frame having conditions over the lattices and accessibility relations.}, and if $w R w'$, then $L_w' \geq^{Cl} L_w$. 
    \item[(2)] $\mathcal{F}_{dec}$: $F \in \mathcal{F}_{dec}$ if, and only if, worlds in $F$ have lattice $T_{(z)}$ for some $z \in \{1, a, b, 0\}$, the accessibility relation $R$ is transitive, and if $w R w'$, then $L_w \geq^{Cl} L_w'$.
\end{enumerate}

The first class of frames have an accessibility relation increasingly classical and the second decreasingly classical. They will validate different formulas in our modal language. For instance, (i) $\mathcal{F}_{inc} \vDash \Box (\varphi \lor \lnot \varphi)$ while (ii) $\mathcal{F}_{dec} \nvDash \Box (\varphi \lor \lnot \varphi)$. To prove (i), we should first note that the value of $\delta \lor \lnot \delta$ is in the filter $Tr$ for all possible valuations of $\delta$ in $T$. Therefore all frames in $\mathcal{F}_{inc}$ or in $\mathcal{F}_{dec}$ validate the formula $\varphi \lor \lnot \varphi$. We further notice that, in a $F \in \mathcal{F}_{inc}$, if $w R w'$, then $L_w' \subseteq L_w$. Hence the value $v_{w'}(\varphi \lor \lnot \varphi)$ will be in the lattice $L_w$ of $w$. It follows that $v_w(\Box (\varphi \lor \lnot \varphi))$ is in the filter -- and $w \vDash \Box (\varphi \lor \lnot \varphi)$. This suffice to show that $\mathcal{F}_{inc} \vDash \Box (\varphi \lor \lnot \varphi)$.

In order to show (ii), we build the counter-model $M$ with two worlds $w$ and $w'$ such that $w R w'$, $L_w = T_{(1)}$, $L_{w'} = T_{(a)}$ and $v_{w'}(\varphi) = \pair{a, a}$:  

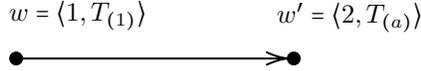
\begin{figure}[H]
\begin{center}
\tikzset{every picture/.style={line width=0.75pt}} 
\begin{tikzpicture}[x=0.75pt,y=0.75pt,yscale=-1,xscale=1]

    \draw  [fill={rgb, 255:red, 0; green, 0; blue, 0 }  ,fill opacity=1 ] (198,140.5) .. controls (198,138.84) and (199.34,137.5) .. (201,137.5) .. controls (202.66,137.5) and (204,138.84) .. (204,140.5) .. controls (204,142.16) and (202.66,143.5) .. (201,143.5) .. controls (199.34,143.5) and (198,142.16) .. (198,140.5) -- cycle ;
    \draw  [fill={rgb, 255:red, 0; green, 0; blue, 0 }  ,fill opacity=1 ] (338,140.5) .. controls (338,138.84) and (339.34,137.5) .. (341,137.5) .. controls (342.66,137.5) and (344,138.84) .. (344,140.5) .. controls (344,142.16) and (342.66,143.5) .. (341,143.5) .. controls (339.34,143.5) and (338,142.16) .. (338,140.5) -- cycle ;
    \draw    (201,140.5) -- (336,140.5) ;
    \draw [shift={(338,140.5)}, rotate = 180] [color={rgb, 255:red, 0; green, 0; blue, 0 }  ][line width=0.75]    (10.93,-3.29) .. controls (6.95,-1.4) and (3.31,-0.3) .. (0,0) .. controls (3.31,0.3) and (6.95,1.4) .. (10.93,3.29)   ;

    \draw (196,109.9) node [anchor=north west][inner sep=0.75pt]    {$w=\langle 1,T_{( 1)} \rangle $};
    \draw (331,110.9) node [anchor=north west][inner sep=0.75pt]    {$w'=\langle 2,T_{( a)} \rangle $};
\end{tikzpicture}
\end{center}
\caption{Counter-model for necessity of excluded middle.}
\label{counter-model-nec-ex-mid}
\end{figure}

In this case $v_w(\Box(\varphi \lor \lnot \varphi)) = \bigcap\limits_{T_{(1)}} \{v_{w'}(\varphi \lor \lnot \varphi)\} = \bigcap\limits_{T_{(1)}} \{\pair{a,a}\}$. Since $\pair{a,a}$ is not in the lattice of $w$, we obtain $v_w(\Box(\varphi \lor \lnot \varphi)) = \bigcap\limits_T \{\pair{a, 1}\} = \pair{a, 1}$. Thence $w \nvDash \Box(\varphi \lor \lnot \varphi)$ as desired.

We further notice that any operation on values of the form $\pair{z,z}$ to themselves results in $\pair{z,z}$. In fact, $\lnot \pair{z,z} = \pair{z,z}$, $\pair{z, z} \lor \pair{z,z} = \pair{z,z}$ and $\pair{z, z} \land \pair{z,z} = \pair{z,z}$. So, for any given formula $\varphi$, there is a valuation of propositional values that results in $\pair{z, z}$ where $z$ is $0$, $a$, $b$ or $1$. To obtain this, one need only to attribute value $\pair{z,z}$ to all propositional variables occurring in $\varphi$. So we can build a counter-example with the frame in \cref{counter-model-nec-ex-mid} for any formula $\Box \varphi$ by ascribing value $\pair{a, a}$ to all propositional variables in $\varphi$. It follows that every formula $\delta$ without modal operators is such that $\mathcal{F}_{dec} \nvDash \Box \delta$. In order to extend this result to formulas in general, one should notice that accessed worlds of accessed worlds will also have the pair $\pair{z,z}$ in their lattice; then one proceeds with a standard proof by induction on formula complexity.

As a general rule, frames in $\mathcal{F}_{dec}$ do not obligate worlds to have necessary formulas. This does not mean that there may not be necessary formulas in particular worlds of a model or necessary formulas in particular worlds. Also, all worlds that only access worlds that are equally classical to themselves will always have many necessary formulas, including $\Box(\varphi \lor \lnot \varphi)$. But if we are `de-classicizing' by going to accessed worlds, it is expected that many sentences will no longer be necessary. The set of possible sentences, nevertheless, is expected to enlarge. Even though we do not have in general $\Box (\varphi \lor \lnot \varphi)$, the decreasing frame validates $\Diamond (\Box (\varphi \lor \lnot \varphi) \lor \lnot \Box (\varphi \lor \lnot \varphi))$ for all worlds that access at least one world (we leave this proof and other characterization results to our ongoing work on the mathematical characterization of these frames). Put another way, if $\mathcal{F}'_{dec}$ is the class of frames in $\mathcal{F}_{dec}$ in which all worlds access at least one world (for example, all worlds access themselves), then $\mathcal{F}'_{dec} \vDash \Diamond (\Box (\varphi \lor \lnot \varphi) \lor \lnot \Box (\varphi \lor \lnot \varphi))$.

Lastly, let us briefly mention some other possible frames one can build in this framework. In the twist structure $T$, we have considered only sublattices that already dismiss `non-contradiction', as pairs $\pair{x, y}$ in frames of $\mathcal{F}_{inc}$ and $\mathcal{F}_{dec}$ do not impose any condition over the product of $x$ and $y$. This is another dimension of becoming more classical. With this in mind, one may generate families of lattices that increase/decrease in classicality both with respect to `excluded-middle' and `non-contradiction'. Moreover, one may consider a dialectic dynamic between increasing and decreasing frames, where if a world accesses two `less-classical' worlds then these two worlds access a common `more-classical' world; for instance:
\begin{enumerate}
    \item[(3)] Dialectic frames $\mathcal{F}_{dial}$: $F \in \mathcal{F}_{dial}$ if, and only if, worlds in $F$ have lattice $T_{(z)}$ for some $z \in \{1, a, b, 0\}$, and
    if $w R w', w''$ and $L_w \geq^{Cl} L_w', L_{w''}$, 
    then there is $w^*$ such that $L_w^* \geq^{Cl} L_w', L_{w''}$ and $w', w'' R w^*$.
\end{enumerate}
Indeed, the addition of a new layer to the concept of frames creates a myriad of possibilities. And in this article we have only considered a very narrow range of the phenomena. We hope, nevertheless, to have provided clear enough picture of how many modal scenarios can be accommodated in many logic modal structures.

\section{Conclusion}

We presented a modal notion and structure that allows putting together modal worlds that operate in different logic systems defined by sublattices of a common lattice. This establishes a natural framework for transferring information between worlds, and addresses the problem of changing meaning explored in the first section. In our system, necessity and possibility of a statement is defined in terms of a comparison between its values in accessible worlds and the common lattice.
We present several examples showing that some properties of necessary$/$possible in standard modal logics are not verified in this broader setting -- e.g. $\varphi$ can be said to be necessary$/$possible even though an/all accessible world falsify $\varphi$.

This setting paves the way for a new field of research. It would be interesting to find conditions over $L$ and its sublattices that ensures traditional properties like the axiom $K$ and others.
We would also like to study strategies to produce the common lattice in our frames from a given set of lattices -- and, consequently, the implications of particular choices on how to put those together. For that purpose, one should look over lattice constructions such as direct product, disjoint union, etc. 
If we allow that our logics can be defined by arbitrary finite linear lattices, then a good candidate to amalgamate any set of these logics is the dense linear order without endpoints given by the limit of Fraïssé's construction of the class of all finite linear orders.

Many variations on the definitions and extensions used in this article should be studied. We have defined only disjunction, conjunction and negation in our frame work, leaving implication defined simply by $-a + b$. One may define (as often done in twist structures) implication as $(-a + b).(--b + -a)$ or else define implication directly using up/down interpretation. Moreover, down-interpretation being a conservative choice with respect to necessity, one may study up-interpretation of necessity as the desired choice for some phenomena.

\printbibliography


\end{document}